\documentclass[12pt,a4paper,reqno]{amsart}
\usepackage{amsmath}
\usepackage{amsfonts}
\usepackage{amssymb}
\numberwithin{equation}{section}

\theoremstyle{plain}

\newtheorem{theorem}[subsection]{Theorem}
\newtheorem{proposition}[subsection]{Proposition}
\newtheorem{lemma}[subsection]{Lemma}
\newtheorem{corollary}[subsection]{Corollary}

\newtheorem{remark}[subsection]{Remark}
\newtheorem{example}[subsection]{Example}

\theoremstyle{definition}

\newtheorem{definition}[subsection]{Definition}

\def\B{{\mathcal B}}

\def\Z{{\mathbb Z}}
\def\E{{\mathbb E}}

\def\R{{\mathbb R}}
\def\N{{\mathbb N}}
\def\Q{{\mathbb Q}}
\def\P{{\mathbb P}}
\newcommand\eps{{\varepsilon}}
\renewcommand\c{{\textbf{c}}}
\newcommand\sgn{{\operatorname{sgn}}}

\begin{document}

\title[Combinatorial and ergodic approaches to Szemer\'edi]{The ergodic and combinatorial approaches to Szemer\'edi's theorem}

\author{Terence Tao}
\address{Department of Mathematics\\University of California at Los Angeles\\ Los Angeles CA 90095}

\email{tao@math.ucla.edu}

\thanks{The author is
supported by a grant from the Packard Foundation.}

\subjclass{11N13, 11B25, 374A5}

\begin{abstract}  A famous theorem of Szemer\'edi asserts that any set of integers of positive upper density will contain arbitrarily long arithmetic progressions.  In its full generality, we know of four types of arguments that can prove this theorem: the original combinatorial (and graph-theoretical) approach of Szemer\'edi, the ergodic theory approach of Furstenberg, the Fourier-analytic approach of Gowers, and the hypergraph approach of Nagle-R\"odl-Schacht-Skokan and Gowers.  In this lecture series we introduce the first, second and fourth approaches, though we will not delve into the full details of any of them.  One of the themes of these lectures is the strong similarity of ideas between these approaches, despite the fact that they initially seem rather different.
\end{abstract}

\maketitle

\section{Introduction}

These lecture notes will be centred upon the following fundamental theorem of Szemer\'edi:

\begin{theorem}[Szemer\'edi's theorem]\cite{szemeredi}  Let $A \subset \Z$ be a subset of the integers of positive upper density, thus $\limsup_{N \to \infty} \frac{|A \cap [-N,N]|}{2N+1} > 0$.  (Here and in the sequel, we use $|B|$ to denote the cardinality of a finite set $B$.)  Then $A$ contains arbitrarily long arithmetic progressions.
\end{theorem}

This theorem is rather striking, because it assumes almost nothing on the given set $A$ - other than that it is large - and concludes that $A$ is necessarily structured in the sense that it contains arithmetic progressions of any given length $k$.
This is a property special to arithmetic progressions (and a few other related patterns).  Consider for instance the question asking whether a set $A$ of positive density must contain a triplet of the form $\{x,y,x+y\}$.  (Compare with the triplet $\{x,y,\frac{x+y}{2}\}$, which is an arithmetic progression of length three.)  It is then clear that the odd numbers, which are certainly a set of positive upper density, do not contain such triples (see however Theorem \ref{schur} below).  Or for another example, consider whether a set of positive upper density must contain a pair $\{x,x+2\}$.  The multiples of $3$ provide an immediate counterexample.  (This is basically why the methods from \cite{gt-primes} can leverage Szemer\'edi's theorem to show that the primes contain arbitrarily long arithmetic progressions, but are currently unable to make any progress whatsoever on the twin prime conjecture.)  But the arithmetic progressions seem to be substantially more ``indestructable'' than these other types of patterns, in that they seem to occur in
any large set $A$ no matter how one tries to rearrange $A$ to eliminate all the progressions.

We have contrasted Szemer\'edi's theorem with some negative results where the selected pattern need not occur.  Now let us 
give the opposite contrast, in which it becomes \emph{very easy} to find a pattern of a certain type in a set.  Here is a basic
example (a special case of a result of Hilbert):

\begin{proposition}\label{xab}  Let $A \subset \Z$ have positive upper density.  Then $A$ contains infinitely many ``parallelograms'' $\{ x, x+a, x+b, x+a+b\}$ where $a,b \neq 0$.
\end{proposition}

Note that if we could just set $a=b$ in these parallelograms then we could find infinitely progressions of length three.  Alas, things are not so easy, and while progressions are certainly intimately related to parallelograms (and more generally to higher-dimensional parallelopipeds, for which an analogue of Proposition \ref{xab} can be easily located), the existence of the latter does not instantly imply the existence of the former without substantial additional effort.  For example, one can easily modify Proposition \ref{xab} to locate, for any $k \geq 1$, infinitely many parallelopipeds of the form $\{ p + \sum_{i \in A} x_i: A \subset \{1,\ldots,k\} \}$ in the primes $\{2,3,5,\ldots\}$, where $p$ is a prime and $x_1,\ldots,x_k > 0$ are positive integers, but this appears to be of no help whatsoever in locating long arithmetic progressions in the primes (one would need to somehow force all the $x_i$ to be equal, which does not seem easily accomplishable).

\begin{proof}  Since $A$ has positive upper density, we can find a $\delta > 0$ and arbitrarily large integers $N$ such that
$$ |A \cap [-N,N]| \geq \delta N.$$
Now consider the collection of all differences $x-y$, where $x,y$ are distinct elements of $A \cap [-N,N]$.  On one hand,
there are $\delta N (\delta N-1)$ possible pairs $(x,y)$ that can generate such a difference.  On the other hand, these differences range from $-2N$ to $2N$, and thus have at most $4N$ possible values.  For $N$ sufficiently large, $\delta N(\delta N-1) > 4N$, and hence by the pigeonhole principle we can find distinct pairs $(x,y), (x',y')$ with $x,y,x',y' \in A \cap [-N,N]$
and $x-y=x'-y'\neq 0$.  This generates a parallelogram.  A simple modification of this argument (which we leave to the reader) in fact generates infinitely many such parallelograms.
\end{proof}

The above argument in fact yields a very large number of parallelograms; if $|A \cap [-N,N]| \geq \delta N$, then $A \cap [-N,N]$ in fact contains $\gg \delta^4 N^3$ parallelograms $\{x,x+a,x+b,x+a+b\}$.  This should be compared against the total number of parallelograms in $[-N,N]$, which is comparable (up to multiplicative constants) to $N^3$.  Thus the density of parallelograms in $A \cap [-N,N]$ differs only by polynomial factors from the density of $A \cap [-N,N]$ itself.  If arithmetic progressions behaved similarly, one would expect a set $A$ in $[-N,N]$ of density $\delta$ to contain $\gg \delta^{C_k} N^2$ arithmetic progressions of a fixed length $k$.  While this is trivially true for $k=2$, it fails even for $k=3$:

\begin{proposition}[Behrend example]\label{behrend}\cite{behrend}  Let $0 < \delta \ll 1$ and $N \geq 1$.  Then there exists a subset $A \subset \{1,\ldots,N\}$ of density $|A|/N \gg \delta$ which contains no more than $\delta^{c \log \frac{1}{\delta}} N^2$ arithmetic
progressions $\{ n,n+r,n+2r \}$ of length three, where $c>0$ is an absolute constant.
\end{proposition}

\begin{proof} The basic idea is to exploit the fact that convex sets in $\R^d$, such as spheres, do not contain arithmetic progressions of length three.  The main challenge is then to somehow ``embed'' $\R^d$ into the interval $\{1,\ldots,N\}$.  To do this, let $M, d \geq 1$ be chosen later, and let $\phi: \{1,\ldots,N\} \to \{0,\ldots,M-1\}^d$ denote the partial base $M$ map $$\phi(n) := ( \lfloor n / M^i \rfloor \operatorname{ mod } M )_{i=0}^{d-1}$$
where $\lfloor x \rfloor$ is the greatest integer less than $x$, and $n \operatorname{ mod } M$ is the remainder of $n$ when divided by $M$.  We then pick an integer $R$ between $1$ and $d M^2$ uniformly at random, and let $B_R \subset \{ 0, \ldots, \lfloor M/10 \rfloor \}^d$ be the set
$$ B_R := \{ (x_1,\ldots,x_d) \in \{ 0, \ldots, \lfloor M/10 \rfloor \}^d: x_1^2 + \ldots + x_d^2 = R \}$$
and then let $A_R := \phi^{-1}(B_R) \subset \{1,\ldots,N\}$ be the preimage of $B_R$.  The set $B_R$ is contained in a sphere and thus contains no arithmetic progressions of length three, other than the trivial ones $\{ x,x,x\}$.  Because there is no ``carrying'' when manipulating base $M$ expansions with digits in $\{ 0, \ldots, \lfloor M/10 \rfloor \}$, we thus conclude
that $A_R$ only contains an arithmetic progression $(n,n+r,n+2r)$ when $r$ is a multiple of $M^d$.  This shows that the number of progressions in $A_R$ is at most $O(M^{-d} N^2)$.  On the other hand, whenever $\phi(n) \in \{0,\ldots,M/10\}^d$, then $n$ has a probability $1/dM^2$ of lying in $A_R$.  Thus we have a lower bound
$$ |A_R| \gg \frac{1}{dM^2} 10^{-d}.$$
If we set $d := c \log \frac{1}{\delta}$ and $M := \delta^c$ for some small constants $c>0$ we obtain the claim.
\end{proof}

This example shows that one cannot hope to prove Szemer\'edi's theorem by an argument as simple as that used to prove Proposition \ref{xab}, as such simple arguments invariably give polynomial type bounds.  Remarkably, this $60$-year old bound of Behrend is still the best known (apart from the issue of optimising the constant $c$).

Another reason why Szemer\'edi's theorem is difficult is that it already implies the much simpler, but still nontrivial,
theorem of van der Waerden:

\begin{theorem}[Van der Waerden's theorem]\label{vdw-thm}\cite{vdw}  Suppose that the integers $\Z$ are partitioned into finitely many colour classes.  Then one of the colour classes contains arbitrarily long arithmetic progressions.
\end{theorem}

Indeed, from the pigeonhole principle one of the colour classes would have positive density, which by Szemer\'edi's theorem gives infinitely long progressions.  The converse deduction is far more difficult; while certain proofs of Szemer\'edi's theorem do indeed use van der Waerden's theorem as a component (e.g. \cite{szemeredi}, \cite{tao:szemeredi}, and Section \ref{szsc} below), many more additional 
arguments are also needed.

While van der Waerden's theorem is not terribly difficult to prove (we give a proof in the next section), it already yields some non-trivial consequences.  Here is one simple one:

\begin{proposition}[Quadratic recurrence]\label{qr}  Let $\alpha$ be a real number and $\eps > 0$.  Then one has $\| \alpha r^2 \|_{\R/\Z} < \eps$ for infinitely many integers $r$, where $\|x\|_{\R/\Z}$ denotes the distance from $x$ to the nearest integer.
\end{proposition}

\begin{proof} Partition the unit circle $\R/\Z$ into finitely many intervals $I$ of diameter $\leq \eps/4$.  Each interval $I$ induces a colour class $\{ n \in \N: \alpha n^2/2 \operatorname{ mod } 1 \in I \}$ on the integers $\Z$.  (This is a basic example of a \emph{structured} colouring; we will see the dichotomy between structure and randomness repeatedly in the sequel.)
By van der Waerden's theorem, one of these classes contains progressions of length $3$ with arbitrarily large spacing $r$, thus for each such $r$ there is an $n$ for which
$$ \alpha n^2/2, \alpha (n+r)^2/2, \alpha(n+2r)^2/2 \in I.$$
The claim now follows from the identity
$$ \alpha n^2/2 - 2\alpha (n+r)^2/2 + \alpha(n+2r)^2/2 = \alpha r^2.$$
\end{proof}

A modification of the argument lets one also handle higher powers $\alpha r^k$.  More general polynomials (with more than one monomial, but with vanishing constant term) can also be handled, although the argument is more difficult.  This simple example already demonstrates however that the number-theoretic question of the distribution of the fractional parts of polynomials is already encoded to some extent within Szemer\'edi's or van der Waerden's theorem.

Szemer\'edi's theorem has many further important extensions and generalisations which we will not discuss here (see for instance Bryna Kra's lectures for some of these). Instead, we will focus on two of the main approaches to proving Szemer\'edi's theorem in its full generality, namely the ergodic theory approach of Furstenberg and the combinatorial approach of R\"odl and coauthors, as well as Gowers.  We will also sketch in very vague terms the original combinatorial approach of Szemer\'edi.  We will however not discuss the important Fourier-analytic approach,
though, despite the many connections between that approach and the ones given here; see Ben Green's lectures for a detailed treatment of the Fourier-analytic method.  
The combinatorial and ergodic approaches
may seem rather different at first glance, but we will try to emphasise the many similarities between them.  In particular, both approaches are based around a \emph{structure theorem}, which asserts that a general object (such as a subset $A$ of the integers) can be somehow split into a ``structured'' component (which has low complexity, is somehow ``compact'', and has high self-correlation) and a ``pseudorandom'' component (which has high complexity, is somehow ``mixing'', and has negligible self-correlation).  One then has to manipulate the structured and pseudorandom components in completely different ways to establish the result.

\section{Prelude: van der Waerden's theorem}

Before we plunge into proofs of Szemer\'edi's theorem, let us first study the much simpler model case of van der Waerden's theorem.  This theorem has both a simple combinatorial proof and a simple dynamical proof; while these proofs do not easily scale up to proving Szemer\'edi's theorem, the comparison between the two is already illustrative.

We begin with the combinatorial proof.  There are three key ideas in the argument (known as a \emph{colour focusing argument}).  The first is to induct on the length of the progression.
The second is to establish an intermediate type of pattern between a progression of length $k$ and a progression of length
$k+1$, which one might call a ``polychromatic fan''.  The third is a concatenation of colours trick in order to leverage the
induction hypothesis on progressions of length $k$, which allows one to move from one fan to the next.

We need some notation.  We use $a+[0,k) \cdot r$ to denote the arithmetic progression $a, a+r, \ldots, a+(k-1)r$.

\begin{definition}  Let $\c: \{1,\ldots,N\} \to \{1,\ldots,m\}$ be a colouring, let $k \geq 1$, $d \geq 0$, and $a \in \{1,\ldots,N\}$.  We define a \emph{fan of radius $k$, degree $d$, and base point $a$} to be a $d$-tuple $(a+[0,k) \cdot r_1, \ldots, a + [0,k) \cdot r_d)$ of progressions in $\{1,\ldots,N\}$ with $r_1,\ldots,r_d>0$.  We
refer to the progressions $a+[1,k) \cdot r_i$, $1 \leq i \leq d$ as the \emph{spokes} of the fan.  We say that a fan
is \emph{polychromatic} if its base point and its $d$ spokes are all monochromatic with distinct colours.   In other words, there exist distinct colours $c_0, c_1, \ldots, c_d \in \{1,\ldots,m\}$ such that
$\c(a) = c_0$, and $\c(a + j r_i) = c_i$ for all $1 \leq i \leq d$ and $1 \leq j \leq k$.
\end{definition}

\begin{theorem}[van der Waerden again]  Let $k, m \geq 1$.  Then there exists $N$ such that any $m$-colouring of
$\{1,\ldots,N\}$ contains a monochromatic progression of length $k$.
\end{theorem}

It is clear that this implies Theorem \ref{vdw-thm}; the converse implication can also be obtained by a simple compactness argument which we leave as an exercise to the reader.

\begin{proof}
We induct on $k$.  The base case $k=1$ is trivial, so suppose $k \geq 2$
and the claim has already been proven for $k-1$.

We now claim inductively that for all $d \geq 0$ there exists a positive integer $N$ such that
any $m$-colouring of $\{1,\ldots,N\}$ contains either a monochromatic progression of length $k$, or
a polychromatic fan of radius $k$ and degree $d$.  The base case $d=0$ is trivial; as soon as
we prove the claim for $d=m$ we are done, as it is impossible in an $m$-colouring for a polychromatic fan
to have degree larger than or equal to $m$.

Assume now that $d > 1$ and the claim has already been proven for $d-1$.  We define $N = 4kN_1N_2$, where $N_1$ and $N_2$ are sufficiently large and will be chosen later.  Let $\c:\{1,\ldots,N\} \to \{1,\ldots,m\}$ be an $m$-colouring of
$\{1,\ldots,N\}$.  Then for any $b \in \{1,\ldots,N_2\}$, the set $\{ bkN_1 + 1, \ldots, bkN_1 + N_1\}$
is a subset of $\{1,\ldots,N\}$ of cardinality $N_1$.  Applying the inductive hypothesis, we see (if $N_1$ is large enough) that
$\{ bkN_1 + 1, \ldots, bkN_1 + N_1\}$ contains either a monochromatic progression of length $k$, or a polychromatic
fan of radius $k$ and degree $d-1$.  If there is at least one $b$ in which the former case applies, we are 
done, so suppose that the latter case applies for every $b$.  This implies that for every
$b \in \{1,\ldots,N_2\}$ there exist $a(b), r_1(b), \ldots, r_{d-1}(b) \in \{1,\ldots,N_1\}$
and distinct colours $c_0(b), \ldots, c_{d-1}(b) \in \{1,\ldots,m\}$ such that $\c(bkN_1+a(b)) = c_0(b)$
and $\c(bkN_1+a(b)+jr_i(b)) = c_i(b)$ for all $1 \leq j \leq k-1$ and $1 \leq i \leq d-1$.  In particular
the map $b \mapsto (a(b), r_1(b), \ldots, r_{d-1}(b),c_0(b), \ldots, c_{d-1}(b))$ is a colouring
of $\{1,\ldots,N_2\}$ by $m^d N_1^d$ colours (which we may enumerate as $\{1,\ldots,m^d N_1^d\}$ in some arbitrary
fashion).  Thus (if $N_2$ is large enough) there exists a monochromatic
arithmetic progression $b+[0,k-1) \cdot s$ of length $k-1$ in $\{1,\ldots,N_2\}$, with some colour
$(a,r_1,\ldots,r_{d-1},c_0,\ldots,c_{d-1})$.  We may assume without loss of generality
that $s$ is negative since we can simply reverse the progression if $s$ is positive.  

Now we use an algebraic
trick (similar to Cantor's famous diagonalization trick) which will convert a progression of identical fans into a new fan of one higher degree, the base points of the original fans being used to form the additional spoke of the new fan.  
Introduce the base point $b_0 := (b-s)kN_1+a$, which lies in $\{1,\ldots,N\}$ by construction of $N$, and consider the fan
$$ (b_0 + [0,k) \cdot skN_1, b_0 + [0,k) \cdot (skN_1 + r_1), \ldots, b_0 + [0,k) \cdot (skN_1 + r_{d-1}))$$
of radius $k$, degree $d$, and base point $b_0$.  We observe that all the spokes of this fan are monochromatic.  For
the first spoke this is because 
$$\c( b_0 + j sk N_1 ) = \c( (b + (j-1) s)kN_1 + a ) = c_0(b + (j-1)s) = c_0$$
for all $1 \leq j \leq k-1$
and for the remaining spokes this is because
$$\c( b_0 + j (sk N_1 + r_t) ) = \c( (b + (j-1) s)kN_1 + a + j r_t ) = c_t(b + (j-1)s) = c_t $$
for all $1 \leq j \leq k-1, 1 \leq t \leq d-1$.
If the base point $b_0$ has the same colour as one of the spokes, then we have found a monochromatic progression of length
$k$; if the base point $b_0$ has distinct colour to all of the spokes, we have found a polychromatic fan of radius $k$
and degree $d$.  In either case we have verified the inductive claim, and the proof is complete.
\end{proof}

Now let us give the dynamical proof.  Van der Waerden's theorem follows from the following abstract topological statement.
Define a \emph{topological dynamical system} to be a pair $(X,T)$ where $X$ is a compact non-empty topological space and $T: X \to X$ is a homeomorphism\footnote{As it turns out, $T$ only needs to be a continuous map rather than a homeomorphism, but we retain the homeomorphism property for some minor technical simplifications.  It is also common to require $X$ to be a metric space rather than a topological one but this does not make a major difference in the argument.}.

\begin{theorem}[Topological multiple recurrence theorem]\label{topmult}\cite{FW}  Let $(X,T)$ be a topological dynamical system.  Then for any open cover $(V_\alpha)_{\alpha \in A}$ of $X$ and $k \geq 2$, at least one of the sets in the cover contains a subset of the form $T^{[0,k) \cdot r} x := \{ x, T^r x, \ldots, T^{(k-1)r} x\}$ for some $x \in X$ and $r > 0$.  (We shall refer to such sets as \emph{progressions of length $k$}.)
\end{theorem}

\begin{proof}[Proof of van der Waerden assuming Theorem \ref{topmult}]  Let $\c: \Z \to \{1,\ldots,m\}$ be an $m$-colouring of the integers.  We can identify $\c$ with a point $x_\c := (\c(n))_{n \in \Z}$ in the discrete infinite product space $\{1,\ldots,m\}^\Z$.  Since each $\{1,\ldots,m\}$ is a compact topological space with the discrete topology, so is $\{1,\ldots,m\}^\Z$.  The shift operator $T: \{1,\ldots,m\}^\Z \to \{1,\ldots,m\}^\Z$ defined by $T( (x_n)_{n \in \Z} ) := (x_{n-1})_{n \in \Z}$ is a homeomorphism.  Let $X$ be the closure of the orbit $\{ T^n x_\c: n \in \Z \}$, then $X$ is also compact, and is invariant under $T$, thus $(X,T)$ is a topological dynamical system.  We cover $X$ by the open sets $V_i := \{ (x_n)_{n \in \Z}: x_0 = i \}$ for $i=1,\ldots,m$; by Theorem \ref{topmult}, one of these open sets, say $V_i$, contains a subset of the form $T^{[0,k) \cdot r} x$
for some $x \in X$ and $r > 0$.  Since $X$ is the closure of the orbit $\{T^n x_\c: n \in \Z\}$, we see from the open-ness of $V_i$ and the continuity of $T$ that $V_i$ must in fact contain a set of the form $T^{[0,k) \cdot r} T^n x_\c$.  But this implies that the progression $-n-[0,k) \cdot r$ is monochromatic with colour $i$, and the claim follows.
\end{proof}

Conversely, it is not difficult to deduce Theorem \ref{topmult} from van der Waerden's theorem, so the two are totally equivalent.  One can view this equivalence as an instance of a \emph{correspondence principle} between colouring theorems and topological dynamics theorems.  By invoking this correspondence principle one leaves the realm of number theory and enters the  infinitary realm of abstract topology.  However, a key advantage of doing this is that we can now manipulate a new object, namely the compact topological space $X$.  Indeed, the proof proceeds by first proving the claim for a particularly simple class of such $X$, the \emph{minimal} spaces $X$, and then extending to general $X$.  This strategy can of course also be applied directly
on the integers, without appeal to the correspondence principle, but it becomes somewhat less intuitive when doing so (we invite the reader to try it!).

The space $X$ encodes in some sense all the ``finite complexity, translation-invariant'' information that is contained in the colouring $\c$.  For instance, if $\c$ is such that one never sees a red integer immediately after a blue integer, this fact will be picked up in $X$ (which will be disjoint from the set $\{ (x_n)_{n \in \Z}: x_0 \hbox{ blue}, x_1 \hbox{ red} \}$).  The correspondence principle asserts that a colouring theorem can be derived purely by exploiting such information. 

\begin{definition}[Minimal topological dynamical system] A topological dynamical system $(X,T)$ is said to be \emph{minimal} if it does not contain any proper subsystem, i.e. there does not exist $\emptyset \subsetneq Y \subsetneq X$ which is closed with $TY=Y$.
\end{definition}

\begin{example} Consider the torus $X = \R/\Z$ with the doubling map $Tx := 2x$.  Then the torus is not minimal, but it contains the minimal system $\{0\}$, the minimal system $\{1/3,2/3\}$, and many other minimal systems.  On the other hand, the same torus with an irrational shift $Tx := x+\alpha$ for $\alpha \not\in \Q$ is minimal.  Minimality can be viewed as somewhat analogous to ergodicity in measure-preserving dynamical systems.  
\end{example}

\begin{lemma} Every topological dynamical system contains at least one minimal topological dynamical subsystem.
\end{lemma}

\begin{proof} Observe that the intersection of any totally ordered chain of topological dynamical systems is again a topological dynamical system (the non-emptiness of such an intersection follows from the finite intersection property of compact spaces).  The claim now follows from Zorn's lemma.
\end{proof}

In light of this lemma, we see that in order to prove Theorem \ref{topmult} it suffices to do so for minimal systems.  One advantage of working with minimal systems is the following.

\begin{lemma}\label{vx} Let $(X,T)$ be a minimal dynamical system, and let $V$ be a non-empty open subset in $X$.  Then $X$ can be covered by finitely many shifts $T^n V$ of $V$.
\end{lemma}

\begin{proof}  If the shifts $T^n V$ do not cover $X$, then the complement $X \backslash \bigcup_{n \in \Z} T^n V$ is a proper closed invariant subset of $X$, contradicting minimality.  Thus the $T^n V$ cover $X$, and the claim follows from compactness.
\end{proof}

\begin{remark} There is a notion of a \emph{minimal colouring} of the integers that corresponds to a minimal system; informally speaking, a minimal colouring is one that does not ``strictly contain'' any other colouring, in the sense that the set of finite blocks of the latter colouring is a proper subset of the set of finite blocks of the former colouring.  This lemma then asserts that in a minimal colouring, any block that does appear in that colouring, in fact appears \emph{syndetically} (the gaps between each appearance are bounded).  Minimal colourings may be considered ``maximally structured'', in that all the finite blocks that appear in the sequence, appear for a ``good reason''.  The opposite extreme to minimal colourings are \emph{pseudorandom colourings}, in which every finite block of colours appears at least once in the sequence (so $X$ is all of $\{1,\ldots,k\}^\Z$).
\end{remark}

Now we can prove Theorem \ref{topmult} for minimal dynamical systems.  We induct on $k$.  The $k=1$ case is trivial; now suppose that $k \geq 2$ and the claim has already been proven for $k-1$, thus given any open cover of $X$, one of the open sets contains a progression of length $k-1$.  Combining this with Lemma \ref{vx} (and the trivial observation that the shift of a progression is again a progression), we obtain

\begin{corollary}\label{xt} Let $(X,T)$ be a minimal dynamical system, and let $V$ be a non-empty open subset in $X$.  Then $V$ contains a progression of length $k-1$.
\end{corollary}

Now we can build fans again.

\begin{definition}  Let $(X,T)$ be a minimal dynamical system, let $(V_\alpha)_{\alpha \in A}$ be an open cover of $X$, let $d \geq 0$, and $x \in X$.  We define a \emph{fan of radius $k$, degree $d$, and base point $x$} to be a $d$-tuple 
$(T^{[0,k) \cdot r_1} x, \ldots, T^{[0,k) \cdot r_d} x)$ of progressions of length $k$ with $r_1,\ldots,r_d > 0$, and
refer to the progressions $a+[1,k) \cdot r_i$, $1 \leq i \leq d$ as the \emph{spokes} of the fan.  We say that a fan
is \emph{polychromatic} if its base point and its $d$ spokes each lie in a distinct element of the cover.  In other words, there exist distinct $\alpha_0,\ldots,\alpha_d \in A$ such that
$x \in A_{\alpha_0}$ and $T^{jr_i} x \in A_{\alpha_i}$ for all $1 \leq i \leq d$ and $1 \leq j \leq k$.
\end{definition}

To prove Theorem \ref{topmult} it now suffices to show

\begin{proposition} Let $(X,T)$ be a minimal dynamical system, and let $(V_\alpha)_{\alpha \in A}$ be an open cover of $X$.
Then for any $d \geq 0$ either there exists at least one polychromatic fan of radius $k$ and degree $d$, or at least one of the sets in the open cover contains a progression of length $k$.
\end{proposition}

Indeed, by compactness we can make the open cover finite, and the above proposition leads to the desired result by taking $d$ large enough.

\begin{proof} The base case $d=0$ is trivial.  Assume now that $d \geq 1$ and the claim has already been proven for $d-1$.  
If one of the $V_\alpha$ contains a progression of length $k$ we are done, so we may assume that we have found a polychromatic fan
$(T^{[0,k) \cdot r_1} x, \ldots, T^{[0,k) \cdot r_{d-1}} x)$ of degree $d-1$, thus 
there exist distinct $\alpha_0,\ldots,\alpha_{d-1} \in A$ such that
$x \in A_{\alpha_0}$ and $T^{jr_i} x \in A_{\alpha_i}$ for all $1 \leq i \leq d-1$ and $1 \leq j \leq k$.
Since the $T^{jr_i}$ are continuous, we can thus find a neighbourhood $V$ of $x$ in $A_{\alpha_0}$ such that
$T^{jr_i} V \subset A_{\alpha_i}$ for all $1 \leq i \leq {d-1}$ and $1 \leq j \leq k$.  By Corollary \ref{xt}
$V$ contains a progression of length $k-1$, say $T^{[1,k) \cdot r_0} y$.  Thus we see that
$T^{jr_0} y \in A_{\alpha_0}$ for $1 \leq j \leq k$, and
$T^{j(r_0+r_i)} y \in A_{\alpha_i}$ for $1 \leq j \leq k$ and $1 \leq i \leq d-1$.  The point $y$ itself lies in
an open set $A_\alpha$.  If $\alpha$ equals one of the $\alpha_0,\alpha_1,\ldots,\alpha_d$, then $V_\alpha$ contains a progression
of length $k$; if $\alpha$ is distinct from $\alpha_0,\alpha_1,\ldots,\alpha_{d-1}$, we have a polychromatic fan of degree $d$.
The claim follows.
\end{proof}

As one can see, the topological dynamics proof contains the same core arithmetical ideas as the combinatorial proof (namely, that a progression of fans can be converted to either a longer progression, or a fan of one higher degree) but the argument is somewhat cleaner as one does not have to keep track of superfluous parameters such as $N$.  For the particular purpose of proving van der Waerden's theorem, the additional overhead in the dynamical proof makes the total argument longer than the combinatorial proof, but for more complicated colouring theorems the dynamical proofs tend to eventually be somewhat shorter and conceptually
clearer than the combinatorial proofs, which often burdened with substantial notation.  The dynamical proofs seem to rely quite
heavily on infinitary tools such as Tychonoff's theorem and Zorn's lemma, though one can reduce the dependence on these tools by making the argument more ``quantitative'' (of course, if one removes the infinitary framework completely, one ultimately ends up at an argument which is more or less just some reworking of the combinatorial argument).

\section{Shelah's argument}

Let us now present another proof of van der Waerden's theorem, due to Shelah \cite{shelah}; it gives slightly better bounds by avoiding inductive arguments which massively increase the number of colours in play.  This argument in fact proves a much stronger theorem, 
namely the \emph{Hales-Jewett theorem}, but we shall content ourselves with a slightly less general result in order to avoid a certain amount
of notation.

\begin{definition}[Cubes]  A \emph{cube} of dimension $d$ and length $k$ is any set of integers of the form
$$ a + [0,k)^d \cdot v = \{ a + n_1 v_1 + \ldots + n_d v_d: 0 \leq n_1,\ldots,n_d\leq k\}$$
where $a \in \Z$ and $v = (v_1,\ldots,v_d)$ is a $d$-tuple of positive integers, with the property that
all the elements $a + n_1 v_1 + \ldots + n_d v_d$ are distinct.
\end{definition}

Cubes are a special case of \emph{generalised arithmetic progressions}, which play an important role in this subject.  

\begin{theorem}[Hales-Jewett theorem]\cite{hales-jewett}\label{hjt}  Let $Q$ be a cube of dimension $d$ and length $k$ which is coloured into $m$ colour classes.  If $j \geq 1$, and
$d$ is sufficiently large depending on $k,m,j$, then $Q$ contains a monochromatic subcube $Q'$ of dimension $j$ and length $k$.
\end{theorem}

Note that the interval $\{1,\ldots,k^d\}$ can be viewed as a proper cube of dimension $d$ and length $k$.  As such, we see that the van der Waerden
theorem follows from the $j=1$ case of the Hales-Jewett theorem.  (The original proof of this theorem proceeded by a colour focusing argument
that directly generalised that used to prove van der Waerden's theorem, and we leave it as an exercise.)

Shelah's proof of this theorem proceeds by an induction on the length $k$.  The $k=1$ case is trivial, so suppose that $k \geq 1$ and that
the theorem has already been proven for $k-1$.  Let us call a subcube
\begin{equation}\label{qp}
 Q' = \{ a + n_1 v_1 + \ldots + n_d v_d: 0 \leq n_1,\ldots,n_d\leq k\}
 \end{equation}
of $Q$ \emph{weakly monochromatic} if whenever one of the $n_1,\ldots,n_d$ is swapped from $k-1$ to $k$ or vice versa, the colour of the element of $Q'$ is unchanged.  It will suffice to show

\begin{theorem}[Hales-Jewett theorem, first inductive step]\label{hjt1}  Let $Q$ be a cube of dimension $d$ and length $k$ which is coloured into $m$ colour classes.  If $j \geq 1$, and $d$ is sufficiently large depending on $k,m,j$, then $Q$ contains a weakly monochromatic subcube $Q'$ of dimension $j$ and length $k$.
\end{theorem}

To prove Theorem \ref{hjt}, one may first without loss of generality ``stretch'' the cube $Q$ by making each $v_i$ enormously large compared
with the previous $v_{i-1}$.  This allows us to eliminate certain ``exotic'' sub-cubes which would cause some technicalities later on.  Then, we let $J$ be a large integer depending on $k,m,j$ to be chosen later.  If $d$ is large enough depending on $k,m,J$, then by Theorem \ref{hjt1} we can find a weakly monochromatic subcube $Q'$ of $Q$ of dimension $J$ and length $k$.  We contract each of the edges by $1$ (deleting all the vertices where one of the $n_i$ is equal to $k$) to create a subcube $Q''$ of $Q$ of dimension $J$ and length $k-1$.  By the induction hypothesis, we see that if $J$ is large enough then $Q''$ will in turn contain a monochromatic cube $Q'''$ of dimension $j$ and length $k-1$.  Since $Q'$ was weakly monochromatic, one can verify that $Q'''$ extends back to a monochromatic cube $Q''''$ of
dimension $j$ and length $k$, which is contained in $Q$, and the claim follows.

It remains to prove Theorem \ref{hjt1}.  Let us modify the notion of weakly monochromatic somewhat.  Let us call the subcube \eqref{qp}
\emph{$i$-weakly monochromatic} for some $0 \leq i \leq d$ if whenever one of the $n_1,\ldots,n_i$ is swapped from $k-1$ to $k$ or vice versa, the colour of the 
element of $Q'$ is unchanged.  It will suffice to show

\begin{theorem}[Hales-Jewett theorem, second inductive step]\label{hjt2}  Let $Q$ be a cube of dimension $d$ and length $k$ which is coloured into $m$ colour classes which is already $i$-weakly monochromatic for some $i \geq 0$.  If $j \geq i+1$, and $d$ is sufficiently large depending on $k,m,j,i$, then $Q$ contains a $i+1$-weakly monochromatic subcube $Q'$ of dimension $j$ and length $k$.
\end{theorem}

Indeed, by iterating Theorem \ref{hjt2} in $i$ we see that for $d$ large enough depending on $k,m,j,i$, $Q$ will contain an $i$-weakly monochromatic subcube of dimension $j$ and length $k$ (the case $i=0$ is trivial); setting $i=j$ we obtain Theorem \ref{hjt1}.

It remains to prove Theorem \ref{hjt2}. As a warmup (and because we need the result to prove the general case) 
let us first give a simple special case of this theorem.

\begin{lemma}[Hales-Jewett theorem, trivial case]\label{hjt3}  Let $Q$ be a cube of dimension $d$ and length $k$ which is coloured into $m$ colour classes.  If $d \geq m+1$, then $Q$ contains a $1$-weakly monochromatic subcube $Q'$ of dimension $1$ and length $k$.
\end{lemma}

\begin{proof}  Write
$$ Q = \{ a + n_1 v_1 + \ldots + n_d v_d: 0 \leq n_1,\ldots,n_d\leq k\}$$
and consider the $m+1$ elements of $Q$ of the form
$$ a + (k-1) v_1 + \ldots + (k-1) v_s + k v_{s+1} + \ldots + k v_{m+1}$$
where $s$ ranges from $1$ to $m+1$.  By the pigeonhole principle two of these have the same colour, thus we have $1 \leq s < s' \leq m+1$
such that the ($1$-dimensional, length $k$) subcube
$$ \{ a + (k-1) v_1 + \ldots + (k-1) v_s + n(v_{s+1} + \ldots + v_{s'}) + k v_{s'+1} \ldots + k v_{m+1}: 1 \leq n \leq k \}$$
is $1$-weakly monochromatic, and the claim follows.
\end{proof}

Now we can prove Theorem \ref{hjt2} and hence the Hales-Jewett theorem.  The main idea is to recast the cube $Q$, not as an $i$-weakly monochromatic $m$-coloured cube of dimension $d$ and length $k$, but rather as an $m^{k^{j-1}}$-coloured cube of dimension $d-j+1$ and length $k$.  More precisely,
let us write
$$ Q = \{ a + n_1 v_1 + \ldots + n_d v_d: 0 \leq n_1,\ldots,n_d\leq k\}$$
and consider now the modified cube of dimension $d-j+1$ and length $k$
$$ \tilde Q := \{ a + n_{j} v_{j} + \ldots + n_d v_d: 0 \leq n_{j},\ldots,n_d\leq k\}.$$
Note that each element $x \in \tilde Q$ is associated to $k^{j-1}$ elements of $Q$, namely
$$ \{ x + n_1 v_1 + \ldots n_{j-1} v_{j-1} \}.$$
Each of these elements has $m$ colours, and so we can naturally associate an $m^{k^{j-1}}$-colouring of $\tilde Q$.  If $d$ (and hence $d-j+1$) is large enough, we can apply Theorem \ref{hjt3} and find a $1$-weakly monochromatic subcube $\tilde Q'$ of dimension $1$ and length $k$ in $\tilde Q$.  It is easy to verify that this in turn induces a $i+1$-weakly monochromatic subcube $Q'$ of dimension $j$ and length $k$ in $Q$, and we are done.

\section{The Furstenberg correspondence principle}

In a previous section, we saw how van der Waerden's theorem was shown to be equivalent to a recurrence theorem in topological dynamics.  Similarly, Szemer\'edi's theorem is equivalent to a recurrence theorem in measure-preserving dynamics.  

\begin{definition}  A \emph{measure-preserving system} $(X, \B, \mu, T)$, is a probability space $(X, \B, \mu)$, where $\B$ is a $\sigma$-algebra of events on $X$, $\mu: \B \to [0,1]$ is a probability measure (thus $\mu$ is countably additive with $\mu(X)=1$), and the \emph{shift map} $T: X \to X$ is a bijection which is bi-measurable (thus $T^n: \B \to \B$ for all $n \in \Z$) and probability preserving (thus $\mu(T^n E) = \mu(E)$ for all $E \in \B$ and $n \in \Z$).
\end{definition}

\begin{example}[Circle shift]  Take $X$ to be the circle $\R/\Z$ with the Borel $\sigma$-algebra $\B$, the uniform probability measure $\mu$, and the shift $T: x \mapsto x+\alpha$ where $\alpha \in \R$.  Thus $T^n E = E+n\alpha$ for any $E \in \B$.
This system is to recurrence theorems as \emph{quasiperiodic sets}, such as the \emph{Bohr set} $\{ n \in \Z: \|n\alpha\|_{\R/\Z} \leq \theta \}$, is to Szemer\'edi's theorem - it is an extreme example of a \emph{structured set}.
\end{example}

\begin{example}[Finite systems]  Take $X$ to be a finite set, and let $\B$ be the $\sigma$-algebra generated by some partition $X = A_1 \cup \ldots \cup A_n$ of $X$ into non-empty sets $A_1,\ldots,A_n$ (these sets are known as ``atoms'').  Thus a set is measurable in $\B$ if and only if it is the finite union of atoms.  We take $\mu$ to be the uniform measure, thus $\mu(E) := |E|/|X|$ for all $E \in \B$.  The shift map $T: X \to X$ is then a permutation on $X$, with the property that it maps atoms to atoms.  Note that if two atoms have different sizes, it will be impossible for the shift map (or any power of the shift map) to take one to the other.  If one assumes that the shift map is \emph{ergodic} (we will define this later), this forces all the atoms to have the same size.  The finite case is not the case of interest in recurrence theorems, but it does serve as a useful toy model that illustrates many of the basic concepts in the proofs without many of the technicalities.  Finite systems have a counterpart in Szemer\'edi's theorem as \emph{periodic sets} - which are trivial for the purpose of demonstrating existence of arithmetic progressions, but still serve as an important illustrative special case for certain components of the proof of Szemer\'edi's theorem.
\end{example}

\begin{remark} The shift $T$ induces an action $n \mapsto T^n$ of the additive integer group $\Z$ on $X$.  One can also study actions of other groups; for instance, actions of $\Z^2$ are described by a pair $S,T$ of commuting bi-measurable probability preserving transformations.
\end{remark}

Given any measure-preserving system $(X,\B,\mu,T)$, a set $E$, and a point $x \in X$, we can define the recurrence set $A = A_{x,E} \subset \Z$ of integers by the formula 
\begin{equation}\label{axe}
A_{x,E} := \{ n \in \Z: T^n x \in E \}.
\end{equation}
This is a way of identifying sets $E$ in a system with sets $A$ in the integers.  Similarly, given a function $f: X \to \R$ on the system, and an $x \in X$, we can define an associated sequence $F = F_{x,f} : \Z \to \R$ by the formula
\begin{equation}\label{fxe}
 F_{x,f}(n) := f(T^n x).
 \end{equation}
This correspondence between sets and functions on the system, and sets and functions on the integers, underlies the Furstenberg correspondence principle.  In particular, it allows one to equate Szemer\'edi's theorem - which is a theorem on the integers - to the following theorem on measure-preserving systems.

\begin{theorem}[Furstenberg multiple recurrence theorem]\label{furst}\cite{furst}  Let $(X,\B,\mu,T)$ be a measure-preserving system.  Then for any set $E \in \B$ of positive measure $\mu(E) > 0$ and any $k \geq 1$, we have
$$ \liminf_{N \to \infty} \E_{1 \leq r \leq N} \mu( E \cap T^r E \cap \ldots \cap T^{(k-1)r} E ) > 0$$
where we use the averaging notation $\E_{1 \leq r \leq N} f(r) := \frac{1}{N} \sum_{r=1}^N f(r)$.
\end{theorem}

\begin{remark} The $k=1$ case is trivial.  The $k=2$ case follows easily from the pigeonhole principle and is known as the \emph{Poincar\'e recurrence theorem}.  The $k=3$ case can be handled by spectral theory (i.e. Fourier analysis).  However the general $k$ case is significantly harder.  It is known that the limit on on the left actually exists, but this is significantly harder (see Bryna Kra's lectures).
\end{remark}

As one consequence of this theorem, we see that every set in $X$ of positive measure contains arbitrarily long progressions.  This should be contrasted with Theorem \ref{topmult}, which can easily be shown to be a special case of Theorem \ref{furst}.

The \emph{Furstenberg correspondence principle} asserts an equivalence between results such as Szemer\'edi's theorem in combinatorial number theory, and recurrence theorems in ergodic theory.  
Let us first show how the recurrence theorem implies Szemer\'edi's theorem.

\begin{proof}[Proof of Szemer\'edi's theorem assuming Theorem \ref{furst}]  This shall be analogous to the topological correspondence principle, in which we shifted the colouring function $\c$ around and took closures to create the dynamical system $X \subset \{1,\ldots,m\}^\Z$.  
This time we shift a set $A$ around and take weak limits to create the measure-preserving system $X \subset \{0,1\}^\Z$.
One can view this as ``inverting'' the correspondence \eqref{axe}; whereas \eqref{axe} starts with a set in a system and turns it into a set of integers, here we need to do things the other way around.

More precisely, suppose for contradiction that Szemer\'edi's theorem fails.  Then there exists a $k \geq 1$, a set $A \subset \Z$ without progressions of length $k$, and a sequence $N_i$ of integers going to infinity such that $\liminf_{i \to \infty} \frac{|A \cap [-N_i,N_i]|}{2N_i+1} > 0$.  Now for each $i$, consider the random set
$$ A_i := A + x_i$$
where $x_i$ is an integer chosen at random from $[-N_i,N_i]$.  As the subsets of $\Z$ can be identified with elements of $X := \{0,1\}^\Z$, we can think of $A_i$ as a random variable taking values in $X$.  More precisely, if we let $\B$ be the Borel $\sigma$-algebra of $X$, we can identify $A_i$ with a probability measure $\mu_i$ on $X$ (it is the average of $2N_i+1$ Dirac masses).  Now $X$ is a separable compact Hausdorff space, and so the probability measures are weakly sequentially compact.  This means that (after passing to a subsequence of $i$ if necessary), the $\mu_i$ converge to another probability measure $\mu$ in the weak sense, thus
$$ \lim_{i \to \infty} \int_X f\ d\mu_i = \int_X f\ d\mu$$
for any continuous function $f$ on $X$.  In particular, if we let\footnote{This is the correct choice of $E$ if one wants to invert the equivalence \eqref{axe}.  Indeed, identifying $A$ with a point in $X$, we see that $A = E_{A,E}$, $A_i = E_{A_i,E}$, and so forth.} $E := \{ (x_n)_{n \in \Z} \in \{0,1\}^\Z: x_n = 1\}$, then
since $E$ is both open and closed,
$$ \lim_{i \to \infty} \mu_i(E) = \mu(E).$$
But a computation shows
$$ \mu_i(E) = \frac{|A \cap [-N_i,N_i]|}{2N_i+1}$$
and hence $\mu(E) > 0$.  Similarly, if $T: X \to X$ is the shift operator $T (x_n)_{n \in \Z} := (x_{n-1})_{n \in \Z}$, then a brief computation shows that
$$ \lim_{i \to \infty} \mu_i( T E ) - \mu_i( E ) = 0$$
and more generally
$$ \lim_{i \to \infty} \mu_i( T F ) - \mu_i( F ) = 0$$
whenever $F$ is a finite boolean combination of $E$ and its shifts.  This means that
$$ \mu(TF) = \mu(F)$$
for all such $F$, and then by the Kolmogorov extension theorem we see that $\mu$ is in fact shift-invariant.  Finally, since $A$ contains no arithmetic progressions of length $k$, we see that
$$ \mu_i( E \cap T^r E \cap \ldots \cap T^{(k-1)r} E ) = 0$$
for any $r > 0$, and hence on taking limits
$$ \mu( E \cap T^r E \cap \ldots \cap T^{(k-1)r} E ) = 0.$$
These facts together contradict the Furstenberg recurrence theorem, and we are done.
\end{proof}

One can easily show that the Szemer\'edi theorem and the Furstenberg recurrence theorem are equivalent to slightly stronger versions of themselves.  For instance, Furstenberg's multiple recurrence theorem generalises to

\begin{theorem}[Furstenberg multiple recurrence theorem, again]\label{furst2}  Let $(X,\B,\mu,T)$ be a measure-preserving system.  Then for any bounded measurable function $f: X \to [0,1]$ with $\int_X f\ d\mu > 0$ and any $k \geq 1$, we have
\begin{equation}\label{liminf}
 \liminf_{N \to \infty} \E_{1 \leq r \leq N} \int_X f T^r f \ldots T^{(k-1)r}f\ d\mu > 0
\end{equation}
where $T^r f := f \circ T^{-r}$ is the translation of $f$ by $r$.
\end{theorem}

This follows simply because if $\int_X f\ d\mu > 0$, then we have the pointwise bound $f \geq c 1_E$ for some $c>0$ and some set $E$ of positive measure, where $1_E$ is the indicator function of $E$.  In a similar spirit, Szemer\'edi's theorem has the following quantitative formulation:

\begin{theorem}[Szemer\'edi's theorem, again]\label{szt} Let $\Z/N\Z$ is a cyclic group.  Then for any bounded function $f: \Z/N\Z \to [0,1]$ with $\E_{n \in \Z/N\Z} f(n) \geq \delta > 0$ and any $k \geq 1$, we have
$$ \E_{n,r \in \Z/N\Z} f(n) T^r f(n) \ldots T^{(k-1)r} f(n) \geq c(k,\delta) $$
for some $c(k,\delta) > 0$ which is independent of $N$, where $T^r f(n) := f(n-r)$.
\end{theorem}

It is easy to see that Theorem \ref{szt} implies Szemer\'edi's theorem in its original formulation, and it can also be easily used (by using the correspondence \eqref{fxe} between functions and sequences) to prove Theorem \ref{furst2} or Theorem \ref{furst} (in fact it gives a lower bound on \eqref{liminf} which depends only on $k$ and the mean $\int_X f\ d\mu$ of $f$).  The converse implication requires an additional averaging argument is essentially due to Varnavides \cite{varnavides}.  We present it here:

\begin{proof}[Proof of Theorem \ref{szt} assuming Szemer\'edi's theorem]  First we observe that for any $k \geq 1$ and $\delta > 0$ that there exists an $M = M(\delta)$ such that any subset of $[1,M]$ of density at least $\delta$ contains at least one progression of length $k$.  For if this were not the case, then one could find arbitrarily large $N$ and sets $A_M \subset [1,M]$ with $|A_M|\geq \delta M$ which contained no progressions of length $k$.  Taking unions of translates of such sets (with $M$ a rapidly increasing sequence) one can easily find a counterexample to Szemer\'edi's theorem.

Now we prove Theorem \ref{szt}.  It is easy to see that $f(n) \geq \delta/2$ on a set $A \subset \Z/N\Z$ of density at least $\delta/2$.  Thus it will suffice to show that
$$ \E_{n,r \in \Z/N\Z} 1_{n,n+r,\ldots,n+(k-1)r \in A} \gg_{k,\delta} 1.$$
For $N$ small depending on $k,\delta$ this is clear (just from taking the $r=0$ case) so assume $N$ is large.  Let $1 \leq M < N$ be chosen later. It will suffice to show that
$$ \E_{n \in \Z/N\Z} \E_{1 \leq r \leq M} 1_{n,n+\lambda r,\ldots,n+(k-1)\lambda r \in A} \gg_{k,\delta} 1$$
for all $\lambda \in \Z/N\Z$, as the claim then follows by averaging in $\lambda$.  We rewrite this as
$$ \E_{n \in \Z/N\Z} \E_{1 \leq m,r \leq M} 1_{n+m,n+m+\lambda r,\ldots,n+m+(k-1)\lambda r \in A} \gg_{k,\delta} 1$$
On the other hand, we have
$$ \E_{n \in \Z/N\Z} \E_{1 \leq m \leq M} 1_{n+\lambda m \in A} = |A|/N \geq \delta$$
so we have $\E_{1 \leq m \leq M} 1_{n+\lambda m \in A} \geq \delta/2$ for a set of $n$ of density at least $\delta/2$.  For each
such $n$, the set $\{ 1 \leq m \leq M: n+\lambda m \in A \}$ has density at least $\delta/2$, and so if we choose $M = M(\delta/2)$ we have at least one $1 \leq m,r \leq M$ for which $n+m,n+m+\lambda r,\ldots,n+m+(k-1)\lambda r \in A$, and so
$$\E_{1 \leq m,r \leq M} 1_{n+m,n+m+\lambda r,\ldots,n+m+(k-1)\lambda r \in A} \gg_M 1.$$
Since $M$ depends on $k,\delta$, the claim follows.
\end{proof}

\begin{remark} One can also deduce Theorem \ref{szt} directly from Theorem \ref{furst} by modifying the derivation of Szemer\'edi's theorem from Theorem \ref{furst}.  We sketch the ideas briefly here.  One can replace $f$ by a set $A$ in $\Z/N\Z$.  One then randomly translates and dilates the function $A$ on $\Z/N\Z$ and then lifts up to $\Z$ to create a random set $A$ in $\Z$.  Now one argues as before.  See \cite{tao:correspondence} for a detailed argument.  See also \cite{bhm} for further exploration of uniform lower bounds in the Furstenberg recurrence theorem.
\end{remark}

\section{Some ergodic theory}

We will not prove Theorem \ref{furst} or Theorem \ref{furst2} here; see Bryna Kra's lectures for a detailed treatment of 
this theory.  However we can illustrate some of the key concepts here.  For those readers which are more comfortable with finite mathematical structures, a good model of a measure-preserving system to keep in mind here is that of the cyclic shift, where $X = \Z/N\Z$, $\B = 2^X$ is the power set of $X$ (so the atoms are just singleton sets) and $T: n \mapsto n+1$ is the standard shift.  
Other finite systems of course exist (though any such system is ultimately equivalent to the disjoint union of finitely many such cyclic shifts).  

The basic ergodic theory strategy in proving Theorem \ref{furst2} is to first prove this result for very \emph{structured} types of functions - functions which have a lot of self-correlation between their shifts.  As it turns out, this is equivalent to studying very structured \emph{factors} $\B'$ of the $\sigma$-algebra $\B$.  One then extends the recurrence result from simple factors to more complicated extensions of these factors, continuing in this process (using Zorn's lemma if necessary) until the full $\sigma$-algebra is recovered (and so all functions are treated).  This is a more complicated version of the topological dynamical situation, in which there was only one type of structured system, namely a minimal system, and the extension from minimal systems to arbitrary systems was trivial (after using Zorn's lemma).

In addition to structured functions, there will also be ``anti-structured'' or ``mixing'' functions which can be considered orthogonal to the structured functions.  These can be viewed as functions for which there is absolutely no correlation between certain of their shifts.  To oversimplify dramatically, one could make the following vague definitions for any $k \geq 2$:

\begin{itemize}
\item A function $f$ is \emph{mixing of order $k-2$} if there is no correlation between the shifts $f$, $T^n f, \ldots, T^{(k-1)n} f$ for generic $n$.
\item A (possibly vector-valued) function $f$ is \emph{strongly structured of order $k-2$} if knowledge of $f$, $T^n f, \ldots, T^{(k-2)n} f$ can be used to predict $T^{(k-1)n} f$ perfectly and ``continuously''.  
\item A function $f$ is \emph{structured of order $k-2$} if it is a component of a strongly structured function of order $k-2$, or can be approximated to arbitrary accuracy by finite linear combinations of such components.
\end{itemize}

These definitions can be formalised, for instance using the Gowers-Host-Kra seminorms; see the lectures of Ben Green and Bryna Kra.  We will not do so here.  However we shall gradually develop some key examples of these concepts in this section.  A fundamental observation in the subject is that there is a \emph{structure theorem} that (for any $k \geq 2$) decomposes any function uniquely into a structured component of order $k-2$ and a mixing component of order $k-2$; indeed, the structured components end up being precisely those functions which are measurable with respect to a special factor $Y_{k-2}$ of $\B$, known
as the \emph{characteristic factor} for $k$-term recurrence\footnote{We are oversimplifying a lot here, there are some subtleties in precisely how to define this factor; in particular the factor $Z_{k-2}$ constructed by Host and Kra \cite{host-kra2} differs slightly from a similar factor $Y_{k-2}$ constructed by Ziegler \cite{ziegler} because a slightly different (but closely related) type of averaging is considered, using $k-1$-dimensional cubes instead of length $k$ progressions.  See  \cite{leibman} for a comparison of the two factors.}.  To prove the Furstenberg recurrence theorem, one first proves recurrence for structured functions of order $d$ for any $d$ (by induction on $d$), and then shows weakly mixing functions of order $k-2$ are negligible for the purpose of establishing $k$-term recurrence.  Setting $d=k-2$ and applying the structure theorem, one obtains the general case.

These matters will be treated in more detail in Bryna Kra's lectures.  Here we shall give only some extremely simple special
cases, to build up some intuition.  There will be a distinct lack of rigour in this section; for instance, we shall omit certain proofs, and be cavalier about whether a function is bounded or merely square integrable, whether a limit actually exists, etc.

We now consider various classes of functions $f: X \to \R$; occasionally we will take $f$ to be complex-valued or vector-valued instead of real-valued.  All functions shall be bounded.

The most structured type of functions $f$ are the \emph{invariant} functions, for which $Tf = f$ (up to sets of measure zero, of course).  These can be viewed as ``(strongly) structured functions of order $0$''.  It is trivial to verify the Furstenberg recurrence theorem for such functions.  It is also clear that these (bounded) functions $f$ form a von Neumann algebra\footnote{It seems clear that the theory of von Neumann algebras is somehow lurking in the background of all of this theory, though strangely enough it does not play a prominent role in the current results. An interesting question is to investigate to what extent this theory would survive if $L^\infty(X)$ was replaced by a \emph{noncommutative} von Neumann algebra.}, as the space $L^\infty(X)^T$ of bounded invariant functions is closed under uniform limits and algebraic operations.  Because of this, we can associate a \emph{factor} $Y_0$ to these functions, defined as the least $\sigma$-algebra with respect to which all functions in $L^\infty(X)^T$ are measurable; because $L^\infty(X)^T$ was a von Neumann algebra, we see that $L^\infty(X)^T$ is in fact \emph{precisely} those functions which are $Y_0$-measurable.  In other words, we take level sets $f^{-1}([a,b])$ of invariant functions and use this to generate the $\sigma$-algebra.  One can equivalently write $Y_0$ as the space of essentially invariant sets $E$, thus $TE$ is equal to $E$ outside of a set of measure zero.  For instance, in the finite case $Y_0$ consists of all sets that are unions of orbits of $T$; in the cyclic case $X = \Z/N\Z$, $Tx = x + n$, $Y_0$ consists of all sets that are cosets of the subgroup generated by $n$ (so if $n$ is coprime to $N$, the only sets in $Y_0$ are the empty set and the whole set).  In the case of the circle shift $X = \R/\Z$, $Tx = x+\alpha$, $Y_0$ is trivial when $\alpha$ is irrational but contains proper subsets of $\R/\Z$ when $\alpha$ is rational.

Complementary to the invariant functions are the \emph{anti-invariant} functions, which are orthogonal to all invariant functions; these are the ``mixing functions of order $0$''.  For instance, given any $g \in L^\infty(X)$, the function $Tg - g$ is an anti-invariant function.  In fact, all invariant functions can be approximated to arbitrary accuracy in $L^2(X)$ as linear combinations of such basic anti-invariant functions $Tg - g$.  This is because if this were not the case, then by the Hahn-Banach theorem there would exist a non-invariant function $f$ which was orthogonal to all of the $Tg - g$.  But then $f$ would be orthogonal to $Tf-f$, which after some manipulation implies that $Tf-f$ has $L^2$ norm zero and so $f$ is invariant, contradiction.  Because of this fact, we see that
anti-invariant functions go to zero in the $L^2$ sense:
\begin{equation}\label{fxt}
f \perp L^\infty(X)^T \implies  \E_{1 \leq r \leq N} T^r f \to_{L^2(X)} 0 \hbox{ as } N \to \infty.
\end{equation}
This can be seen by first testing on basic anti-invariant functions $Tg-g$ (in which case one has a telescoping sum), 
taking linear combinations, and then taking limits.  One specific consequence of this is the mixing property
\begin{equation}\label{gvn-0}
\lim_{N \to \infty} \E_{1 \leq r\leq N} \int_X f T^r g\ d\mu = 0
\end{equation}
whenever \emph{at least one} of $f$ and $g$ is anti-invariant.  (Note that there is a symmetry due to the identity
$\int_X f T^r g = \int_X g T^{-r} f$.)  We will refer to this as the \emph{generalised von Neumann theorem of order $0$}.

From Hilbert space theory we know that every function $f$ in $L^2(X)$ uniquely splits as the sum of an invariant function $f_{U^\perp}$ and an anti-invariant $f_{U}$ function.  In fact, since the invariant functions are not only a closed subspace of $L^2(X)$, but are also the measurable functions with respect to a factor $Y_0$, we can write explicitly $f_{U^\perp} = \E(f|Y_0)$ and $f_U = f - \E(f|Y_0)$, where the \emph{conditional expectation operator} $f \mapsto \E(f|Y_0)$ is simply the orthogonal projection from $L^2(X)$ to the subspace $L^2(Y_0)$ of $Y_0$-measurable functions.

If $f$ is invariant, then clearly its averages converge back to $f$:
$$ f \in L^\infty(X)^T \implies  \E_{1 \leq r \leq N} T^r f \to_{L^2(X)} f \hbox{ as } N \to \infty.$$
Combining this with \eqref{fxt} (and taking limits to extend $L^\infty$ to $L^2$)
we obtain the \emph{von Neumann ergodic theorem}
$$ f \in L^2(X) \implies \E_{1 \leq r \leq N} T^r f \to_{L^2(X)} \E(f|Y_0) \hbox{ as } N \to \infty.$$
This implies in particular that
$$ \E_{1 \leq n \leq N} \int_X f T^n f\ d\mu \to \int_X f \E(f|Y_0)\ d\mu = \| \E(f|Y_0) \|_{L^2(Y_0)}^2$$
which already proves the $k=2$ case of the Furstenberg recurrence theorem (and gives a precise value for the limit).

\begin{example} Consider the case of finite systems.  Then the invariant functions are those functions which are constant of each of the orbits of $T$, while the anti-invariant functions are those functions which have mean zero on each of the orbits of $T$.  If $f: X \to \R$ is a general function, then the invariant part $\E(f|Y_0)$ is the function which assigns to each orbit of $T$ (i.e. to each atom of $Y_0$) the average value of $f$ on that orbit, while the anti-invariant part $f - \E(f|Y_0)$ is formed by subtracting the mean of each orbit from the original function.  It is an instructive exercise to verify all the arguments used to prove the von Neumann ergodic theorem directly in this finite system case.
\end{example}

The factor $Y_0$ also leads to a useful \emph{ergodic decomposition} of a general measure-preserving system into ergodic ones.  A measure preserving system is said to be \emph{ergodic} if $Y_0$ is trivial, thus every invariant set has measure either zero or one (or equivalently that every invariant function is constant almost everywhere).  One can view the space $X$ and the $\sigma$-algebra $\B$ as fixed, in which case ergodicity is a property of the shift-invariant probability measure $\mu$.  Then it turns out that while a general measure $\mu$ is not ergodic, it can always be decomposed (or \emph{disintegrated}) as an integral $\int_Y \mu_y\ d\nu(y)$ of ergodic shift-invariant probability measures $\mu_x$ parameterised by some parameter $y$ on another probability space $Y$.  To formalise this decomposition in general requires a certain amount of measure theory, but in the case of a finite system the process is quite simple to describe.  Namely, take $Y$ to be the system $(X,Y_0,\mu)$, and for each $y \in Y$ let $\mu_y$ be the uniform distribution on the $T$-orbit $\{ T^n y: n \in \Z\}$ of $y$.  Then one easily verifies that $\mu = \int_Y \mu_y\ d\nu(y)$, and that each $\mu_y$ is an ergodic measure (all invariant sets either have zero measure or full measure).
The ergodic decomposition in this case is essentially just the decomposition of $X$ into individual orbits of $T$, upon each of which $T$ is ergodic.  One can easily use the ergodic decomposition to reduce the task of proving Furstenberg's recurrence theorem to the special case in which the system is ergodic; we omit the details.  This is somewhat analogous to the reduction in topological dynamics to minimal systems.  Unfortunately, whereas in the dynamical case the assumption of minimality was very strong and lead quickly to a proof of the topological recurrence theorem, ergodicity is not by itself a strong enough condition to quickly obtain a direct proof of Furstenberg's recurrence theorem, and further classification and decomposition of the measure-preserving system is needed.  As it turns out, one usually cannot usefully disintegrate the measure $\mu$ into any smaller invariant measures once one is at an ergodic system; however it is still possible (and useful) to disintegrate the measures into \emph{non-invariant} measures, where the shift map does not act separately on each component, but instead mixes them together using something called a ``cocycle''.  A simple finitary example occurs when considering a finite ergodic system $(X,\B,\mu,T)$ with $\B = 2^X$ which contains a shift-invariant factor $\B' \subset \B$.  The ergodicity forces all the atoms in $\B'$ to be the same size, and thus they are all bijective (non-canonically) to a single set $Z$.  This allows one can then parameterise $X$ as $Y \times Z$, where $Y$ is the collection of all the atoms of $\B'$; since the shift $T$ maps one such atom to another, the factor $(X,\B',\mu,T)$ is then equivalent to a system $(Y, 2^Y, \nu, S)$ on $Y$ where $\nu$ is uniform measure on $Y$, and the original shift can then be described as $T(y,z) := (Sy, \rho_y(z))$ where for each $y \in Y$, the \emph{cocycle} $\rho_y: Z \to Z$ is a permutation on $Z$.  One can view $X$ as an \emph{extension} of $Y$, by converting each point $y$ to a ``vertical fiber'' $y \times Z$.  We can disintegrate $\mu = \int_Y \mu_y\ d\nu(y)$ where $\mu_y$ is uniform measure on $\{y\} \times Z$.  These measures are not invariant; instead $T$ will map $\mu_y$ to $\mu_{Sy}$ for all $y$. The iterates $T^n$ are then described as $T^n(y,z) = (S^n y, \rho_{y,n}(z))$, where the $\rho_{y,n}$ are defined using the \emph{cocycle equation}
$$ \rho_{y,n+m} = \rho_{S^m y, n} \circ \rho_{y,m}.$$
This is a more complicated version of the more familiar equation $T^{n+m} = T^n \circ T^m$, thus cocycles are more complicated versions of shifts (indeed as we just saw, a cocycle is simply the ``vertical component'' of a shift in a larger product space).  The study of cocycles forms an integral part of the higher order recurrence theory but will not be discussed here.

Now let us look at double recurrence (the $k=3$ case of Theorem \ref{furst2}), in which we investigate the limiting behavior 
of averages such as
\begin{equation}\label{ftf}
\lim_{N \to \infty} \E_{1 \leq r \leq N} \int_X f T^r f T^{2r} f\ d\mu.
\end{equation}
If $f$ is invariant, then again this expression is easy to compute (it is just $\int_X f^3$).  One may hope, as in the preceding discussion, that anti-invariant functions are negligible, in the sense that
$$ \lim_{N \to \infty} \E_{1 \leq r \leq N} \int_X f T^r g T^{2r} h\ d\mu = 0$$
whenever $f,g,h$ are bounded at least one of $f,g,h$ is anti-invariant.  Unfortunately, this is not the case.  For a very simple example, take the small cyclic group $X = \Z/M\Z$ for odd $M$ and let $f=g=h$ be the function which equals $M-1$ at $0$ and $-1$ elsewhere.  Then these functions are all anti-invariant, but the above average can be computed to be $M^2 - 1$; the problem is that periodically (whenever $n$ is a multiple of $M$) there is a huge ``spike'' in the value of $\int_X f T^n g T^{2n} h\ d\mu$ which imbalances the average dramatically.  Thus periodic functions (ones in which $T^n f = f$ for some $n > 0$) cause a problem.  More generally\footnote{A simple application of Fourier analysis or the spectral theorem reveals that every periodic function is a finite linear combination of eigenfunctions, with eigenvalues equal to roots of unity.}, the \emph{eigenfunctions}, in which $T f = e^{2\pi i\theta} f$ for some $\theta\in \R/\Z$, will also cause a problem (note that invariant functions correspond to the case $\theta=0$).  Indeed if one sets $h:=f$ and $g := \overline{f}^2$, then we
see that $T^{2r} h = e^{4\pi ir\theta} h$ and $T^r g = e^{-4\pi ir\theta} g$, and hence\footnote{This corresponds to the fact that sets of integers such as the \emph{Bohr set} $\{ n \in \Z: \| \alpha n \|_{\R/\Z} \leq \eps\}$ have an unexpectedly high number of progressions of length three, due to the identity $\alpha n - 2 \alpha(n+r) +\alpha(n+2r) = 0$, which implies that if two elements of a progression lie in the Bohr set, then the third element has an unexpectedly high probability of doing so also.  One should caution that this is not always the case; with the Behrend example in Proposition \ref{behrend}, when two elements of a progression lie in the set, then the third element has an unexpectedly \emph{small} probability of lying in the set.  Thus certain types of structure can in fact \emph{reduce} the number of progressions present, though Szemer\'edi or Furstenberg tells us that they cannot destroy these progressions completely.  This is another indication that the proof of this theorem has to be somewhat nontrivial (in particular, a naive symmetrisation or variational argument will not work).}
$$  \lim_{N \to \infty} \E_{1 \leq r \leq N} \int_X f T^r f T^{2r} f\ d\mu = \int_X |f|^4\ d\mu \neq 0,$$
despite the fact that such eigenfunctions will necessarily be anti-invariant for $\theta \neq 0$ (as eigenfunctions of the unitary operator $T$ with distinct eigenvalues are necessarily orthogonal).  

However, one can simply deal with these problems by devising a suitable factor (larger than $Z_0$) to contain them.  For instance, one can create the factor $Y_0$ generated by all the periodic functions.  This factor can be larger than $Z_0$ (e.g. in the finite case, $Y_0$ is in fact everything).  The periodic functions form an algebra (they are closed under arithmetic operations) but are not quite a von Neumann algebra because they are not quite closed under limits\footnote{There does not seem to be a conventional name for what the uniform or $L^2$ limit of periodic functions should be called.  One possibility is ``pro-periodic'' or ``profinitely periodic'' functions.}.  Nevertheless, the periodic functions are still \emph{dense} in $L^2(Z_0)$, which turns out to be good enough for most purposes.  Even larger than $Y_0$ is $Z_1$, the factor generated by all eigenfunctions - this factor is known as the \emph{Kronecker factor}.  Now the eigenfunctions are not closed under addition (though they are closed under multiplication), however the space of \emph{quasiperiodic functions} - finite linear combinations of eigenfunctions - is indeed an algebra.  The closure of the quasiperiodic functions in $L^2$ are the \emph{almost periodic functions} - and this is a von Neumann algebra, indeed an $L^2$ function is almost periodic if it is measurable in $Z_1$.  One can classify all these properties in terms of the orbit $\{ T^n f: n \in \Z\}$:
\begin{itemize}
\item $f$ is invariant if and only if the orbit $\{ T^n f: n \in \Z\}$ is a singleton.
\item $f$ is periodic if and only if the orbit $\{ T^n f: n \in \Z\}$ is finite.
\item $f$ is an eigenfunction if and only if the orbit $\{ T^n f: n \in \Z\}$ lives in a one-dimensional complex vector space.
\item $f$ is quasiperiodic if and only if the orbit $\{ T^n f: n \in \Z\}$ lives in a finite-dimensional vector space.
\item $f$ is almost periodic if and only if the orbit $\{ T^n f: n \in \Z\}$ is precompact (its closure is compact).
\end{itemize}

Functions in these classes will be referred to as ``structured functions of order $1$'' or ``linearly structured functions''; the eigenfunctions\footnote{An individual quasiperiodic function is usually not strongly structured, in the sense that $f(x)$ does not determine $T^n f(x)$ in a continuous manner; however a quasiperiodic function is the component of a vector-valued function which is strongly structured.  For instance, if $X = (\R/\Z)^2$ and $T(x_1,x_2) = (x_1+\alpha_1,x_2+\alpha_2)$ for rationally independent $\alpha_1,\alpha_2$, then $f(x_1,x_2) := e^{2\pi i (x_1 + x_2})$ is quasiperiodic but not strongly structured, however the vector-valued function $(e^{2\pi i(x_1+x_2)}, e^{2\pi i x_1}, e^{2\pi i x_2})$ \emph{is} strongly structured.} are ``strongly structured functions of order $1$''.  
The linear comes from the fact that the action of $T^n f$ behaves ``linearly'' in $n$; observe for instance that if $f$ is an eigenfunction with eigenvalue $e^{2\pi i\theta}$ then $T^n f = e^{2\pi i n \theta} f$.
Now it turns out that one can get a good handle on the average \eqref{ftf} for all $f$ in the linearly structured classes - and more precisely we have a non-trivial lower bound when $f$ is non-negative and not identically zero.  We already saw what happened when $f$ was invariant.  If instead $f$ was periodic with some period $m$, then we get a large positive contribution to \eqref{ftf} (specifically, $\int_X f^3\ d\mu$) when $n$ is a multiple of $m$, which is already enough for a non-trivial lower bound.  For the other cases, one can use a pigeonhole argument to show that almost periodic functions behave very much like periodic functions (hence the name), in the sense that given any $\eps$, we have $\|T^n f - f \|_{L^2(X)} \leq\eps$ for a set of $n$ of positive density.  Note that if $T^n f$ is close to $f$, then (by applying $T^n$ and then the triangle inequality) $T^{2n} f$ is close to $f$ also, which can be used (together with H\"older's inequality and the boundedness of $f$) to show that $f T^n f T^{2n} f$ is close to $f^3$.  This gives a contribution close to $\int_X f^3\ d\mu$ for all $n$ in a set of positive density, and one still gets a good lower bound for $f$.  Note that these arguments extend easily to higher averages such as those involving
$\int_X f T^n f \ldots T^{(k-1)n} f\ d\mu$.  (But problems will emerge with the other half of the argument, as orthogonality to linear structure is not enough to eliminate all problems with triple and higher recurrence.)

There is another proof of recurrence for almost periodic functions which looks more complicated, but ends up being more robust and can extend (with some effort) to higher order cases.  We know that the orbit $\{ T^n f: n \in \Z\}$ is precompact, which means that for any $\eps > 0$ one can cover this orbit by finitely many balls.  This allows us to apply the van der Waerden theorem (or its topological counterpart) and conclude the existence of many progressions $n, n+r, \ldots, n+(k-1)r$ for which $T^n f, T^{n+r} f, \ldots, T^{n+(k-1)r} f$ are all close to each other.  This means that $\int_X f T^r f \ldots T^{(k-1)r} f\ d\mu$ is close to $\int_X f^k\ d\mu > 0$, which can be used as before to get a nontrivial lower bound.

Now we say that a function $f$ is ``mixing of order $1$'', or ``linearly mixing'', 
if it is orthogonal to all almost periodic functions, or in other words $\E(f|Z_1) = 0$.  It turns out that a more useful characterisation of this mixing property exists.

\begin{lemma}\label{wmix}  A real-valued function $f \in L^\infty(X)$ is mixing of order $1$ if and only if the self-correlation functions $T^n f f$ are asymptotically mixing of order $0$, in the sense that
\begin{equation}\label{limn}
 \lim_{N \to \infty} \E_{-N \leq n \leq N} \| \E(T^n f f | Z_0) \|_{L^2}^2 = 0.
\end{equation}
\end{lemma}

\begin{proof} (Sketch only)  Suppose first that $f$ obeys the property \eqref{limn}. A Cauchy-Schwarz argument (based on something called the \emph{van der Corput lemma}), which we omit, then shows that
$$  \lim_{N \to \infty} \E_{-N \leq n \leq N} \| \E(T^n g f | Z_0) \|_{L^2}^2 = 0$$
for any bounded $g$.  If we apply this in the particular case that $g$ is an eigenfunction, we have
$\| \E(T^n g f | Z_0) \|_{L^2} = \| \E(gf|Z_0) \|_{L^2}$ and hence $\E(gf|Z_0) = 0$ for all eigenfunctions $g$.  In particular $f$ is orthogonal to all eigenfunctions, hence to all quasiperiodic functions, hence to all almost periodic functions, and is thus mixing of order $1$.

Now suppose that \eqref{limn} fails.  We rewrite the left-hand side (ignoring issues regarding interchange of limit and integral, which can be justified using the von Neumann ergodic theorem applied to the product space $X \times X$) as
$$ \langle f, \lim_{N \to\infty} \E_{-N \leq n \leq N} \E(T^n f f|Z_0) T^n f \rangle.$$
Let us introduce the linear operator $S: L^2(X) \to L^2(X)$ by
$$ S f := \lim_{N \to\infty} \E_{-N \leq n \leq N} \E(T^n f g|Z_0) T^n f$$
(again, let us ignore the issue regarding whether this limit exists).  Thus $\langle f, Sf \rangle \neq 0$.
This is a self-adjoint operator (in fact, it is positive definite).  Also, being the limit of averages of finite rank operators, it can be shown to be a compact operator.  Finally, we have the translation invariance property $T^n S = S T^n$.  In particular, this shows that the orbit of $Sf$ lies in the range of $S$ and is thus precompact:
$$ \{ T^n Sf: n \in \Z \} = \{ S T^n f : n \in \Z \} \subset \{ S g: \|g\|_{L^2(X)} \leq\|f\|_{L^2(X)} \}.$$
This shows that $Sf$ is almost periodic.  Thus $f$ is not orthogonal to all almost periodic functions, a contradiction.
\end{proof}

By using \eqref{limn} and some Cauchy-Schwarz (more precisely, using the van der Corput lemma) one can show that weakly mixing functions of order $1$ are negligible for the purposes of double recurrence; indeed, we have
$$ \lim_{N \to \infty} \E_{1 \leq r \leq N} \int_X f T^r g T^{2r} h\ d\mu = 0$$
whenever $f,g,h$ are bounded and at least one of $f,g,h$ are mixing of order $1$.  We can refer to this as the \emph{generalised von Neumann theorem of order $1$}.  On the other hand, every bounded function $f$ has a unique decomposition $f = \E(f|Z_1) + (f - \E(f|Z_1)$ as an almost periodic function $\E(f|Z_1)$ and a weakly mixing function $f - \E(f|\Z_1)$; I like to refer to this as the \emph{Koopman-von Neumann theorem}\footnote{Lemma \ref{wmix} is also sometimes known as the Koopman-von Neumann theorem; the two facts are of course closely related.}.  Note also that if $f$ is non-negative with positive mean, then the almost periodic component $\E(f|Z_1)$ will be also.  Combining this fact with the recurrence already obtained for almost periodic functions, and
the negligibility of weakly mixing functions, we obtain recurrence for all functions, i.e. we have established the general $k=3$ case of Furstenberg's multiple recurrence theorem.

We now give the barest sketch of how things continue onward from here.  For $k=4$ one needs to define notions of almost periodicity and weak mixing of order $2$.  Of the two, the latter is easier, because we can copy Lemma \ref{wmix}, and declare a function $f$ to be weakly mixing of order $2$ if its self-correlations $T^n f f$ are asymptotically weakly mixing of order $1$, thus 
$$  \lim_{N \to \infty} \E_{-N \leq n \leq N} \| \E(T^n f f | Z_1) \|_{L^2}^2 = 0.$$
(Many other equivalent definitions are possible.)  Repeated application of van der Corput eventually shows that such functions are negligible for the averages
$$ \lim_{N \to \infty} \E_{1 \leq r \leq N} \int_X f T^r g T^{2r} h T^{3r} k\ d\mu$$
in the sense that this average vanishes whenever $f,g,h,k$ are bounded and at least one is weakly mixing of order $2$.  It is not hard to show that there exists a unique factor $Z_2$ (that extends $Z_1$) such that the weakly mixing functions of order $2$ are precisely those functions $f$ whose conditional expectation $\E(f|Z_2)$ vanishes.  (In the work of Host and Kra, this factor $Z_2$ is generated by \emph{nonconventional averages} such as
$$ \lim_{N \to \infty} \E_{-N\leq a,b,c \leq N} T^a f T^b f T^c f T^{a+b} f T^{b+c} f T^{a+c} f T^{a+b+c} f;$$
this idea was then adapted for the finite setting in \cite{gt-primes} as the notion of a \emph{dual function} to construct a finitary analogue of this factor.)  One would then like the almost periodic functions of order $2$ to be some dense subclass of
$L^2(Z_2)$.  This can be done; the trick is to repeat the original definition of almost periodic, but view terms such as ``finite dimensional'' or ``compact'' not in terms of vector spaces over $\R$ (as we have implicitly been doing), but rather\footnote{The combinatorial analogue of this would be to partition the original space $X$ into atoms - in this case, the atoms of $Z_1$, and somehow work on each atom separately.  Of course, things are not this simple because the atoms are usually not shift-invariant and so the shift structure is now more complicated, passing from one atom to the next.  The graph theoretic approach, which we will discuss later, also relies heavily on restriction to atoms, but can cope with this with much greater ease because this approach ``forgets'' all the arithmetic structure and so there is nothing to destroy when passing to an atom.} as \emph{modules} over the von Neumann algebra $L^\infty(Z_1)$ of bounded almost periodic functions.  In particular:

\begin{itemize}
\item $f$ is an \emph{eigenfunction of order $2$} (also known as a \emph{quadratic eigenfunction}) if and only if the orbit $\{ T^n f: n \in \Z\}$ lives in a one-dimensional module over $L^\infty(Z_1)$.
\item $f$ is \emph{quasiperiodic of order $2$} if and only if the orbit $\{ T^n f: n \in \Z\}$ lives in a finite-dimensional module over $L^\infty(Z_1)$.
\item $f$ is \emph{almost periodic of order $2$} if and only if the orbit $\{ T^n f: n \in \Z\}$ can be ``approximated to arbitrary accuracy'' by subsets of finite-dimensional modules over $L^\infty(Z_1)$.  (The precise definition is a little tricky and subtle; see \cite{furst-book}.)
\end{itemize}

A quadratic eigenfunction can equivalently be defined (at least in the ergodic case) as a function $f$ obeying an identity of the form $Tf = gf$, where $g$ is itself a linear eigenfunction, thus $Tg = e^{2\pi i\theta} g$ for some $\theta \in \R/\Z$.  The origin of the term ``quadratic'' can then be observed from an inspection of the phase in the identity
$$ T^n f = e^{2\pi i n(n-1) \theta} g^n f.$$
From the closely related identity
$$ f T^n (\overline{f}^3) T^{2n} (f^3) T^{3n} \overline{f} = |f|^8$$
one also sees that quadratic eigenfunctions are not negligible for the purposes of triple recurrence (indeed they end up being orthogonal to all quadratically mixing functions). Quasiperiodic functions of order $2$ are special cases of \emph{$2$-step nilsequences}, which will be discussed in Bryna Kra's lectures.  They can be viewed as components of \emph{vector-valued} (or matrix-valued) quadratic eigenfunctions, and arise from what are known as \emph{finite rank extensions} of the Kronecker factor $Z_1$.

At any rate, the almost periodic functions of order $2$ now form a dense subclass of $L^2(Z_2)$, and are an algebra, and so one can repeat previous arguments and reduce the proof of the Furstenberg recurrence theorem for $k=3$ to the task of proving such recurrence for such quadratically almost periodic functions.  This turns out to be complicated - in part because this result includes Proposition \ref{qr} as a special case (the case of quadratic eigenfunctions), and this proposition is itself not
entirely trivial (requiring at a bare minimum some form of van der Waerden's theorem).  Fortunately, the colouring argument given previously for almost periodic functions - which \emph{does} use van der Waerden's theorem - extends (after nontrivial effort)
to this case, and more generally to all orders, thus leading to a proof of the Furstenberg recurrence theorem.  See \cite{furst}, \cite{FKO}, \cite{furst-book}, as well as Bryna Kra's lectures.

\section{The graph theoretic approach}

Now we leave ergodic theory and turn to what (at first glance) appears to be a completely different approach to Szemer\'edi's theorem, though at a deeper inspection one will find many themes in common.  In the ergodic approach, it was the shift operator $T$ which was the primary focus of investigation; the underlying set $A$ of integers merely provided some probability measure for $T$ to leave invariant.  We have seen that the dynamical approach focuses almost entirely on the shift operator.  In marked contrast,
the hypergraph approach discards the shift structure completely; instead, it views the problem of finding an arithmetic progression as that of solving a set of simultaneous relations; these relations initially have some additive structure, but this structure is soon discarded, as these relations are soon modeled abstractly by graphs and hypergraphs.  With the forgetting of so much structure it is remarkable that any nontrivial progress can still be made; however there turn out to be deep theorems 
in (hyper)graph theory, comparable (though not directly equivalent) to the deep recurrence theorems in topological dynamics and 
ergodic theory, which allow one to proceed even after losing almost all of the arithmetic structure.  It is a fascinating
question as to what the ``true'' origin of these deep facts are - it seems to be some very abstract and general dichotomy between randomness and structure - and how they may be united with the ergodic and Fourier-analytic approaches.

To illustrate the power of the graph theoretic approach, let us prove a theorem which looks similar to van der Waerden's theorem though it is slightly different.

\begin{theorem}[Schur's theorem]\label{schur}  Suppose the positive integers $\Z_+$ are finitely coloured.  Then one of the colour classes contains a triple of the form $\{x,y,x+y\}$.
\end{theorem}

\begin{proof} Our task is to find $x,y > 0$ and a colour class ${\mathcal C}$ for which we have the simultaneous relations
\begin{align*}
x &\in {\mathcal C}\\
y &\in {\mathcal C}\\
x+y &\in {\mathcal C}.
\end{align*}
The problem is that these equations (three relations in two unknowns) are coupled together in an unpleasant way. However we can decouple things slightly by making the (somewhat underdetermined) substitution $x = b-a$, $y = c-b$ for some $a < b < c$; our task is then to find such $a<b<c$ and a colour class ${\mathcal C}$ for which we have the simultaneous relations
\begin{align*}
b-a &\in {\mathcal C}\\
c-b &\in {\mathcal C}\\
c-a &\in {\mathcal C}.
\end{align*}
Now we have three relations in three unknowns, which is a bit better for the purposes of finding solutions.  Furthermore, the relations are more symmetric in $a,b,c$, and each relation only involves two of the three unknowns.  This is all that we will need to proceed.  Indeed, let us now edge-colour the complete graph on the natural numbers by assigning to each edge $(a,b)$ with $b>a$, the colour of $b-a$ in the original colouring (this is known as the \emph{Cayley graph} associated to the original colouring).  A solution to the above simultaneous relations is now nothing more than a monochromatic triangle in this graph.  But the existence of such a triangle follows immediately from Ramsey's theorem.  (Indeed one sees that one can even take $a,b,c$ to be no larger than $6$!)
\end{proof}

Note that we only used a very special case of Ramsey's theorem; using the full version of Ramsey's theorem leads to substantial generalisation of Schur's theorem, especially when combined with van der Waerden's theorem, known as \emph{Rado's theorem}; see for instance \cite{graham}.

Now we see what can similarly be done for progressions of length three in a set $A$ of integers.  Actually it will be convenient to localise to a cyclic group $\Z/N\Z$ and prove the following.

\begin{theorem}[Roth's theorem, cyclic group version]\label{cyc}  Let $N$ be a large integer, and let $A \subset \Z/N\Z$ be such that $|A| \geq \delta N$. Then there are at least $c(\delta) N^2$ progressions $x,x+r,x+2r$ in $A$ for some $c(\delta) >0$ (we allow $r$ to be zero).
\end{theorem}

It is easy to see that this implies the $k=3$ version of Szemer\'edi's theorem (and is in fact equivalent to it, thanks to the formulation in Theorem \ref{szt}).
Our task is to find many solutions to the system of relations
\begin{align*}
n &\in A\\
n+r &\in A\\
n+2r &\in A
\end{align*}
Again this is three equations in two unknowns.  We add an unknown by making the underdetermined substitution $n := -x_2-2x_3$, $r := x_1+x_2+x_3$ and obtain the system
$$
\begin{array}{llll}
 & -x_2 &-2x_3 &\in A \\
x_1 & & - x_3 &\in A \\
-2x_1 & - x_2 & & \in A
\end{array}
$$
This is again three relations in three unknowns, where each relation involves only two of the three variables; our task is to locate $c(\delta) N^3$ solutions.  The situation is not quite the same as with Schur's theorem, though; for instance, the three relations are not entirely symmetric.  On the other hand, we already know a lot of \emph{degenerate} solutions to this system:
$$
\begin{array}{llll}
 & -x_2 &-2x_3 &\in A \\
x_1 & & - x_3 &\in A \\
-2x_1 & - x_2 & & \in A\\
x_1 & + x_2 & + x_3 & = 0.
\end{array}
$$
Indeed, every element of $A$ generates $N$ such solutions, so we have $\delta N^2$ solutions in all.  We can rephrase this as a conditional probability bound
\begin{equation}\label{xxx}
 \P( -x_2-2x_3, x_1-x_3, 2x_1+x_2 \in A | x_1+x_2+x_3 = 0 ) \geq \delta
 \end{equation}
where we think of $x,y,z$ as ranging freely over the cyclic group $\Z/N\Z$, and then conditioned so that $x+y+z=0$.  Our goal seems innocuous, namely to remove this conditional expectation and conclude that
\begin{equation}\label{xxx-2}
\P( -x_2-2x_3, x_1-x_3, 2x_1+x_2 \in A ) \geq c(\delta).
\end{equation}
This is less trivial than it first appears.  The problem is that the event $x+y+z=0$ has tiny probability - $1/N$ - and so we only get a tiny lower bound of $\delta/N$ if we naively apply Bayes' identity.  (This corresponds to the fact that the number of trivial progressions - $\delta N^2$ - is negligible compared with the number of progressions that we actually want, which is $c(\delta) N^3$.)  However, the point will be that the solution set $\{ (x,y,z): -x_2-2x_3, x_1-x_3, 2x_1+x_2 \in A\}$, being the intersection of three ``second-order'' sets $\{ (x_1,x_2,x_3): -x_2-2x_3 \in A\}$, $\{ (x_1,x_2,x_3): -x_1-x_3 \in A\}$, $\{ (x_1,x_2,x_3): 2x_1+x_2 \in A\}$,
is not a completely arbitrary set, and as it turns out it cannot concentrate itself entirely on the ``third-order set''
$\{ (x_1,x_2,x_3): x_1+x_2+x_3 = 0\}$.  For instance, observe that given \emph{any} relation $x_i \sim x_j$ involving just two of the $x_1,x_2,x_3$, we have
\begin{equation}\label{pxx}
\P( x_i \sim x_j | x_1+x_2+x_3 = 0 ) = \P( x_i \sim x_j )
\end{equation}
or given any sets $A_1,A_2$, we have
$$ \P( x_1 \in A_1, x_2 \in A_2 | x_1+x_2 = 0 ) \leq \min( \P( x_1 \in A_1), \P( x_2 \in A_2) ) 
\leq \P( x_1 \in A_1, x_2 \in A_2 )^{1/2}.$$
So we see that when the structure of the set is sufficiently ``low order'', one can remove the conditional expectation.  Can one do so here?  The answer is yes, and it relies on the following abstract result.

\begin{lemma}[Triangle removal lemma]\label{trl}\cite{rsz} Let $G$ be a graph on $n$ vertices that contains fewer than $\eps n^3$ triangles for some $0 < \eps < 1$.  Then it is possible to delete $o_{\eps \to 0}(n^2)$ edges from $G$ to create a triangle-free graph $G'$.
\end{lemma}

As usual we use $o_{\eps \to 0}(X)$ to denote a quantity which is bounded by $c(\eps) X$ for some function $c(\eps)$ of $\eps$ which goes to zero as $\eps \to 0$.  Later on we will allow the decay rate to depend on additional parameters, for instance $o_{\eps \to 0;k}(1)$ would be a quantity which decayed to zero as $\eps \to 0$ for each fixed $k$, but which need not decay uniformly in $k$.  An equivalent formulation of this lemma is:

\begin{lemma}[Triangle removal lemma, again]\label{trl2} Let $G$ be a graph on $n$ vertices that contains at least $\delta n^2$ edge-disjoint triangles for some $0 < \delta < 1$.  Then it must in fact contain $c(\delta) n^3$ triangles, where $c(\delta) > 0$ depends only on $\delta$.
\end{lemma}

We leave the equivalence of these two formulations to the reader.  From the second formulation it is an easy matter to deduce
\eqref{xxx-2} from \eqref{xxx}, by considering the tripartite graph formed by three copies of $V$ (corresponding to $x_1,x_2,x_3$ respectively), and with the three edge classes between these copies defined by the relations
$-x_2-2x_3 \in A$, $x_1-x_3 \in A$, and $2x_1+x_2 \in A$ respectively; again, we leave this as an exercise for the reader.

There is another way to phrase this lemma in a ``several variable measure theory'' language that brings it more into line with the ergodic theory approach (and also the Fourier-analytic approach).

\begin{lemma}[Triangle removal lemma, several variable version]\label{trl3} Let $(X,\mu_X)$, $(Y,\mu_Y)$, $(Z,\mu_Z)$ be probability spaces, and let $f: X \times Y \to [0,1]$, $g: Y \times Z \to [0,1]$, and $h: Z \times X \to [0,1]$ be measurable functions such that
$$ \Lambda_3(f,g,h) \leq \eps$$
for some $0 < \eps < 1$, where $\Lambda_3$ is the trilinear form
$$ \Lambda_3(f,g,h) := \int_X \int_Y \int_Z f(x,y) g(y,z) h(z,x)\ d\mu_X(x) d\mu_Y(y) d\mu_Z(z).$$
Then there exists functions 
$\tilde f: X \times Y \to [0,1]$, $\tilde g: Y \times Z \to [0,1]$, and $\tilde h: Z \times X \to [0,1]$ which differ
from $f,g,h$ in $L^1$ norm by $o_{\eps \to 0}(1)$, thus
\begin{align*}
 \int_X \int_Y |f(x,y)-\tilde f(x,y)|\ d\mu_X(x) d\mu_Y(y), &
\int_Y \int_Z |g(y,z)-\tilde g(y,z)|\ d\mu_Y(y) d\mu_Z(z),&
\int_Z \int_X |h(z,x)-\tilde h(z,x)|\ d\mu_Z(z) d\mu_X(x)& \leq o_{\eps \to 0}(1),
\end{align*}
and such that $\tilde f(x,y) \tilde g(y,z) \tilde h(z,x)$ vanishes identically (in particular, $\Lambda_3(\tilde f, \tilde g, \tilde h) = 0$).
\end{lemma}

One can easily deduce Lemma \ref{trl} from Lemma \ref{trl3} by specialising $X,Y,Z$ to be the finite vertex set $V$ with the uniform probability measure, and let $f=g=h$ be the indicator function of the edge set of the graph $G$; we omit the details.  The converse implication is also true but somewhat tricky (one must discretise the measure spaces $X,Y,Z$, and split the atoms of such spaces to approximate the probability measures by uniform distributions, and also replace the functions $f,g,h$ by indicator functions); we again omit the details.  We will choose to work with the analytic formulation of the triangle removal lemma in these notes because it seems to extend more easily to the hypergraph setting (in which one considers similar expressions in more variables, where now each function can depend on three or more variables).

Lemma \ref{trl3} asserts, roughly speaking, that if a collection of low complexity functions have a small product, then one can ``clean'' each function slightly \emph{in a low-complexity manner} in order to make the product vanish entirely.  Note that the claim would be trivial if one were allowed to modify (say) $f$ in a manner which could depend on all three variables $x,y,z$.  The power of the lemma lies in the fact that the high-complexity expression $\Lambda_3(f,g,h)$ can be manipulated purely in terms of low-complexity operations.  This rather deep phenomenon seems to be rather general; in fact there is a similar lemma for any non-negative combination of functions of various collections of variables (we shall describe one such version a little later below).  It is however still not perfectly well understood.

The way one proves Lemma \ref{trl3} is by decomposing $f,g,h$ into ``structured'' or ``low complexity'' components, which are easier to clean up, and ``error terms'', which for one reason or another do not interfere with the cleaning process because they give a negligible contribution to expressions such as $\Lambda_3(f,g,h)$.
It turns out that there are two types of error terms which come into play.  The first are errors which are ``small'' in an integral sense, say in $L^2$ norm, while the second are errors which are (very) small in a weak sense (for instance, they are small when tested against other functions which depend on other sets of variables).  The latter will be encoded using a useful norm, the \emph{Gowers $\Box^2$ norm} $\|f\|_{\Box^2(X \times Y)} = \|f\|_{\Box^2}$, defined for measurable bounded $f: X \times Y \to \R$ by the formula
$$ \|f\|_{\Box^2(X \times Y)}^4 := \int_X \int_X \int_Y \int_Y f(x,y) f(x,y') f(x',y) f(x',y')\ d\mu_X(x) d\mu_X(x') d\mu_Y(y) d\mu_Y(y').$$
One easily verifies that the right-hand side is non-negative.  From two applications of the Cauchy-Schwarz inequality one verifies the \emph{Gowers-Cauchy-Schwarz inequality}
\begin{equation}\label{gczi} 
\begin{split}
|\int_X \int_X \int_Y \int_Y &f_{00}(x,y) f_{01}(x,y') f_{10}(x',y) f_{11}(x',y')\ d\mu_X(x) d\mu_X(x') d\mu_Y(y) d\mu_Y(y')|\\
&\leq \|f_{00}\|_{\Box^2} \|f_{01}\|_{\Box^2} \|f_{10}\|_{\Box^2} \|f_{11}\|_{\Box^2}
\end{split}
\end{equation}
from which one readily verifies that $\Box^2$ obeys the triangle inequality and is thus at least a seminorm.  From the Gowers-Cauchy-Schwarz inequality (and bounding the $\Box^2$ norm crudely by the $L^\infty$ norm) one also sees that
\begin{equation}\label{gcz}
|\int_X \int_Y f(x,y) g(y) h(x)\ d\mu_Y(x) d\mu_Y(y)| \leq \|f\|_{\Box^2}
\end{equation}
whenever $g, h$ are measurable functions bounded in magnitude by $1$; this in particular shows that if $\|f\|_{\Box^2}=0$ then $f$ is zero almost everywhere.  Thus the $\Box^2$ norm is indeed a norm\footnote{One can also identify the $\Box^2$ norm with the Schatten-von Neumann $4$-norm of the integral operator with kernel $f(x,y)$; in the important special case when $X$ is a finite set with the uniform distribution, and $f$ is symmetric, then the $\Box^2$ norm is simply the $l^4$ norm of the eigenvalues of the matrix associated to $f$.  If $f$ is the indicator function of a graph $G$, the $\Box^2$ norm is a normalised count of the number of $4$-cycles in $G$.  However we will not take advantage of these facts as they do not generalise well to hypergraph situations.}, after the customary convention of identifying two functions that agree almost everywhere.  Letting $g,h$ depend on a third variable $z$ in \eqref{gcz} and integrating in $z$, and using symmetry, we thus conclude the \emph{generalised von Neumann inequality}
\begin{equation}\label{gvni}
 |\Lambda_3(f,g,h)| \leq \min( \|f\|_{\Box^2}, \|g\|_{\Box^2}, \|h\|_{\Box^2} )
\end{equation}
whenever $f: X \times Y \to [-1,1]$, $g: Y \times Z \to [-1,1]$, $h: Z \times X \to [-1,1]$ are measurable.
 
Thus functions with tiny $\Box^2$ norm have a negligible impact on the $\Lambda_3$ form; such functions are known as \emph{pseudorandom} or \emph{Gowers uniform}.  To exploit this, one would now like
to decompose arbitrary functions $f: X \times Y \to [0,1]$ into a ``structured'' component which can be easily analysed and manipulated, plus errors which are small in $\Box^2$ or are otherwise easy to deal with.  The first key observation is

\begin{lemma}[Lack of uniformity implies correlation with structure]\label{corstruct}  Let $f: X \times Y \to [-1,1]$ be such that $\|f\|_{\Box^2} \geq \eta$ for some $\eta > 0$.  Then there exists $A \subset X$ and $B \subset Y$ such that
$$ |\int_X \int_Y 1_A(x) 1_B(y) f(x,y)\ d\mu_X(x) d\mu_Y(y)| \geq \eta^4/4.$$
\end{lemma}

\begin{proof} By definition of the $\Box^2$ norm we have
$$ \int_X \int_Y \int_X \int_Y f(x,y) f(x,y') f(x',y) f(x',y')\ d\mu_X(x) d\mu_X(x') d\mu_Y(y) d\mu_Y(y') \geq \eta^4.$$
By the pigeonhole principle and the boundedness of $f$, we can thus find $x',y'$ such that
$$ |\int_X \int_Y f(x,y) f(x,y') f(x',y)\ d\mu_X(x) d\mu_Y(y)| \geq \eta^4.$$
We rewrite this using Fubini's theorem as
$$ |\int_{-1}^1 \int_{-1}^1 \sgn(s) \sgn(t) \int_X \int_Y 1_{A_s}(x) 1_{B_t(y)} f(x,y)\ d\mu_X(x) d\mu_Y(y) ds dt| \geq \eta^4$$
where $A_s := \{ x \in X: \sgn(s) f(x,y') \geq |s|\}$ and $B_t := \{ y \in Y: \sgn(t) f(x',y) \geq |t|\}$.  The claim then follows
from another application of the pigeonhole principle.
\end{proof}
 
To exploit this we borrow some notation from the ergodic theory approach, namely that of $\sigma$-algebras and conditional
expectation.  However, in this simple context we will only need to deal with \emph{finite} $\sigma$-algebras.  If $\B$ is a finite factor of $X$ (i.e. a finite $\sigma$-algebra of measurable sets in $X$), then $\B$ is essentially just a partition of $X$ into finitely many disjoint atoms $A_1,\ldots,A_M$ (more precisely, $\B$ is the $\sigma$-algebra consisting of all finite unions of these atoms).  If $f: X \to \R$ is measurable, then the conditional expectation $\E(f|\B): X \to \R$ is the function defined as $\E(f|\B)(x) := \frac{1}{A_i} \int_{A_i} f(x)\ d\mu_X(x)$ whenever $x$ lies in an atom $A_i$ of positive measure.  (Conditional expectations are only defined up to sets of measure zero, so we can define $\E(f|\B)$ arbitrarily on atoms of measure zero.)  We say that a factor has \emph{complexity at most $m$} if it is generated by at most $m$ sets (and thus it contains at most $2^m$ atoms).  If $\B_X$ is a finite factor of $X$ with atoms $A_1,\ldots,A_M$, and $\B_Y$ is a finite factor of $Y$ with atoms $B_1,\ldots,B_N$, then $\B_X \vee \B_Y$ is a finite factor of $X \times Y$ with atoms $A_i \times B_j$ for $1 \leq i \leq M$ and $1 \leq j \leq N$.

The key relationship between the $\Box^2$ norm and conditional expectation on finite factors is the following.

\begin{lemma}[Lack of uniformity implies energy increment]\label{fei}  Let $\B_X, \B_Y$ be finite factors of $X,Y$ respectively of complexity at most $m$, and let $f: X \times Y \to [0,1]$ be such that
$$ \|f - \E(f|\B_X \vee \B_Y) \|_{\Box^2(X \times Y)} \geq \eta$$
for some $\eta > 0$.  Then there exists extensions $\B'_X$, $\B'_Y$ of $\B_X, \B_Y$ of complexity at most $m+1$ such that
$$ \|\E(f|\B'_X \vee \B'_Y) \|_{L^2(X \times Y)}^2 \geq \|\E(f|\B_X \vee \B_Y) \|_{L^2(X \times Y)}^2 + \eta^8/16.$$
Here of course $\|F\|_{L^2(X \times Y)}^2 := \int_X \int_Y |F(x,y)|^2\ d\mu_X(x) d\mu_Y(y)$.
\end{lemma}

The key point here is that $f$ - which is a ``second-order'' object, depending on two variables - is correlating with
two ``first-order'' objects $\B'_X$, $\B'_Y$.  This ultimately will allow us to approximate the second-order object by a number of first-order objects.  It is this kind of reduction - in which a single high-order object is traded in for a large number of lower-order objects - which is the key to proving results such as the triangle removal lemma.  The quantity
$\|\E(f|\B_X \vee \B_Y) \|_{L^2(X \times Y)}^2$ is known as the \emph{index} of the partition $\B_X \vee \B_Y$ in the graph theory literature; here we shall refer to it as the \emph{energy} of this partition.

\begin{proof}  From Lemma \ref{corstruct} we can find measurable $A \subset X$, $B \subset Y$ such that
$$ | \int_X \int_Y 1_A(x) 1_B(y) (f - \E(f|\B_X \vee \B_Y))\ d\mu_X(x) d\mu_Y(y) | \geq \eta^4/4.$$
Let $\B'_X$ be the factor of $X$ generated by $\B_X$ and $A$, and similarly let $\B'_Y$ be the factor of $Y$ generated
by $\B_Y$ and $B$, then $\B'_X, \B'_Y$ have complexity at most $m+1$.  Since $1_A(x) 1_B(y)$ is $\B'_X \vee \B'_Y$ measurable, we have
\begin{align*}
&\int_X \int_Y 1_A(x) 1_B(y) (f - \E(f|\B_X \vee \B_Y))\ d\mu_X(x) d\mu_Y(y) =\\
&\quad \int_X \int_Y 1_A(x) 1_B(y) \E(f - \E(f|\B_X \vee \B_Y) | \B'_X \vee \B'_Y )\ d\mu_X(x) d\mu_Y(y)
\end{align*}
so by Cauchy-Schwarz
$$\| \E(f - \E(f|\B_X \vee \B_Y) | \B'_X \vee \B'_Y ) \|_{L^2(X \times Y)} \geq \eta^4/4.$$
Now observe that the quantity
$$ \E(f - \E(f|\B_X \vee \B_Y) | \B'_X \vee \B'_Y ) = \E(f|\B'_X \vee \B'_Y) - \E(f|\B_X \vee \B_Y) $$
is orthogonal to $\E(f|\B_X \vee \B_Y)$.  The claim then follows from Pythagoras' theorem.
\end{proof}

Note that if $f$ is bounded by $1$, then the quantity $\|\E(f|\B_X \vee \B_Y) \|_{L^2(X \times Y)}^2$ is bounded between $0$ and $1$.  Thus an easy iteration of the above lemma gives

\begin{corollary}[Koopman-von Neumann decomposition]\label{kv}  Let $\B_X, \B_Y$ be finite factors of $X,Y$ respectively of complexity at most $m$, let $f: X \times Y \to [0,1]$ be measurable, and let $\eta > 0$.  Then there exists extensions $\B'_X$, $\B'_Y$ of $\B_X, \B_Y$ of complexity at most $m+\frac{16}{\eta^8}$ such that
$$  \|f - \E(f|\B'_X \vee \B'_Y) \|_{\Box^2(X \times Y)} < \eta.$$
\end{corollary}

This corollary splits $f$ into a bounded complexity object $\E(f|\B'_X \vee \B'_Y)$ and an error which is small in the $\Box^2$ norm.  In practice, this decomposition is not very useful because the complexity of the structured component $\E(f|\B'_X \vee \B'_Y)$ is large compared to the bounds available on the error $f - \E(f|\B'_X \vee \B'_Y)$.  However one can rectify this by one further iteration of the above decomposition:

\begin{lemma}[Szemer\'edi regularity lemma]\label{srl}  Let $f: X \times Y \to [0,1]$ be measurable, let $\tau > 0$, and let $F: \N \to \N$ be an arbitrary increasing function (possibly depending on $\tau$).  Then there exists an integer $M = O_{F,\tau}(1)$ and a decomposition $f = f_1 + f_2 + f_3$ where
\begin{itemize}
\item ($f_1$ is structured) We have $f_1 = \E(f|\B_X \vee \B_Y)$ for some finite factors $\B_X, \B_Y$ of $X,Y$ respectively of complexity at most $M$;
\item ($f_2$ is small) We have $\|f_2\|_{L^2(X \times Y)} \leq \tau$.
\item ($f_3$ is very uniform) We have $\|f_3\|_{\Box^2(X \times Y)} \leq 1/F(M)$.
\item (Positivity) $f_1$ and $f_1+f_2$ take values in $[0,1]$.
\end{itemize}
\end{lemma}

This lemma may not immediately resemble the usual Szemer\'edi regularity lemma for graphs, but it can easily be used to deduce that lemma.  See \cite{tao:revisited}.  One can obtain a result similar to this from spectral theory, by viewing $f$ as the kernel of an integral operator and decomposing $f$ using the singular value decomposition of that operator, with $f_1,f_2,f_3$ corresponding to the high, medium, and low singular values respectively.  However it then takes some effort to ensure that $f_1$ and $f_1+f_2$ are non-negative.  See \cite{greenkon} for some related discussion. The more ``ergodic'' approach here, relying on conditional expectation, gives worse quantitative bounds but does easily ensure the positivity property, which is crucial in many applications.

\begin{proof}  Construct recursively a sequence of integers
$$ 0 = M_0 \leq M_1 \leq M_2 \leq \ldots$$
by setting $M_0 := 0$ and $M_i := M_{i-1} + 16 F(M_{i-1})^8$ for $i \geq 1$.  Then for each $i \geq 0$, construct recursively
factors $\B^i_X, \B^i_Y$ of $X,Y$ of complexity at most $M_i$ by setting $\B^0_X$ and $\B^0_Y$ to be the trivial factors of complexity $0$, and then applying Corollary \ref{kv} repeatedly to let $\B^i_X, \B^i_Y$ be extensions of $\B^{i-1}_X$, $\B^{i-1}_Y$
such that
$$ \|f - \E(f|\B^i_X \vee \B^i_Y) \|_{\Box^2(X \times Y)} < 1/F(M_{i-1}).$$
The energies $\|\E(f|\B^i_X \vee \B^i_Y) \|_{L^2(X \times Y)}^2$ are monotone increasing in $i$ by Pythagoras' theorem, and
are bounded between $0$ and $1$.  Thus by the pigeonhole principle we can find $1 \leq i \leq 1/\tau^2$ for which
$$ \|\E(f|\B^i_X \vee \B^i_Y) \|_{L^2(X \times Y)}^2 \leq \|\E(f|\B^{i-1}_X \vee \B^{i-1}_Y) \|_{L^2(X \times Y)}^2 + \tau^2.$$
If one then sets
\begin{align*}
f_1 &:= \E(f|\B^{i-1}_X \vee \B^{i-1}_Y); \\
f_2 &:= \E(f|\B^i_X \vee \B^i_Y) - \E(f|\B^{i-1}_X \vee \B^{i-1}_Y); \\
f_3 &:= f - \E(f|\B^i_X \vee \B^i_Y); \\
M &:= M_{i-1}
\end{align*}
then we see that the claims are easily verified.
\end{proof}

A slight modification of the above argument allows one to simultaneously regularise several functions at once using the same partition.  More precisely, we have

\begin{lemma}[Simultaneous Szemer\'edi regularity lemma]\label{ssrl}  Let $f: X \times Y \to [0,1]$, $g: Y \times Z \to [0,1]$, $h: Z \times X \to [0,1]$ be measurable, let $\tau > 0$, and let $F: \N \to \N$ be an arbitrary increasing function (possibly depending on $\tau$).  Then there exists an integer $M = O_{F,\tau}(1)$, factors $\B_X, \B_Y, \B_Z$ of $X,Y,Z$ respectively of complexity at most $M$ and decompositions 
$f = f_1 + f_2 + f_3$, $g = g_1 + g_2 + g_3$, $h = h_1 + h_2 + h_3$, where
\begin{itemize}
\item ($f_1$, $g_1$, $h_1$ are structured) We have $f_1 = \E(f|\B_X \vee \B_Y)$, $g_1 = \E(g|\B_Y \vee \B_Z)$, and
$h_1 = \E(f|\B_Z \vee \B_X)$.
\item ($f_2$, $g_2$, $h_2$ are small) We have $\|f_2\|_{L^2(X \times Y)}, \|g_2\|_{L^2(Y \times Z)}, \|h_2\|_{L^2(Z \times X)} \leq \tau$.
\item ($f_3$, $g_3$, $h_3$ are very uniform) We have $\|f_3\|_{\Box^2(X \times Y)}, \|g_3\|_{\Box^2(X \times Y)}, \|h_3\|_{\Box^2(X \times Y)} \leq 1/F(M)$.
\item (Positivity) $f_1,g_1,h_1$ and $f_1+f_2,g_1+g_2,h_1+h_2$ take values in $[0,1]$.
\end{itemize}
\end{lemma}

We leave the proof of this lemma as an exercise to the reader.  With this lemma we can now prove Lemma \ref{trl3}.  Actually we shall prove a slightly stronger statement, which provides more information about the functions $\tilde f$, $\tilde g$, $\tilde h$ involved.

\begin{lemma}[Strong triangle removal lemma, several variable version]\label{trl4} Let $(X,\mu_X)$, $(Y,\mu_Y)$, $(Z,\mu_Z)$ be probability spaces, and let $f: X \times Y \to [0,1]$, $g: Y \times Z \to [0,1]$, and $h: Z \times X \to [0,1]$ be measurable functions such that $\Lambda_3(f,g,h) \leq \eps$ for some $0 < \eps < 1$. Then there exists factors $\B_X, \B_Y, \B_Z$ of $X,Y,Z$ respectively of complexity at most $O_{\eps}(1)$ and sets $E_{X,Y} \in \B_X \vee \B_Y$, $E_{Y,Z} \in \B_Y \vee \B_Z$, $E_{Z,X} \in \B_Z \vee \B_X$ respectively
with $1_{E_{X,Y}}(x,y) 1_{E_{Y,Z}}(y,z) 1_{E_{Z,X}}(z,x)$ vanishing identically, such that
\begin{align*}
\int_X \int_Y f(x,y) 1_{E_{X,Y}^c}(x,y)\ d\mu_X(x) d\mu_Y(y),&\\
\int_Y \int_Z g(y,z) 1_{E_{Y,Z}^c}(y,z)\ d\mu_Y(y) d\mu_Z(z),&\\
\int_Z \int_X h(z,x) 1_{E_{Z,X}^c}(z,x)\ d\mu_Z(z) d\mu_X(x) &\leq o_{\eps \to 0}(1).
\end{align*}
\end{lemma}

Note that Lemma \ref{trl4} immediately implies Lemma \ref{trl3} by setting $\tilde f := f 1_{E_{x,y}}$, etc.  This strengthened version of the lemma will come in handy in the next section.

\begin{proof}  We apply Lemma \ref{ssrl} with $0 < \tau \ll 1$ and $F$ to be chosen later; for now, one should think of $\tau$ as being moderately small, but not very small compared to $\eps$, and similarly $F$ will be a moderately growing function.  This gives us an integer $M = O_{F,\tau}(1)$, factors $\B_X,\B_Y,\B_Z$ of complexity at most $M$, and decompositions $f=f_1+f_2+f_3$, etc. with the stated properties.  In particular
$$ \Lambda_3(f_1+f_2+f_3, g_1+g_2+g_3, h_1+h_2+h_3) \leq \eps.$$
The idea shall be to eliminate the uniform errors $f_3,g_3,h_3$, and then the small errors $f_2,g_2,h_2$, leaving one with only the structured components $f_1,g_1,h_1$, which will be easy to deal with directly.

It is easy to eliminate $f_3,g_3,h_3$.  Indeed from repeated application of the generalised von Neumann inequality \eqref{gvni} and the $\Box^2$ bounds on $f_3,g_3,h_3$ we have
\begin{equation}\label{fgh}
 \Lambda_3(f_1+f_2, g_1+g_2, h_1+h_2) \leq \eps + O( 1 / F(M) ).
\end{equation}
We would now like to similarly eliminate $f_2, g_2, h_2$.  A naive application of the $L^2$ bounds would give an estimate of the form
\begin{equation}\label{fgh-bad}
 \Lambda_3(f_1, g_1, h_1) \leq \eps + O(\tau) + O( 1 / F(M) )
\end{equation}
but the $O(\tau)$ error turns out to be far too expensive for our purposes.  Instead we proceed in a more ``local'' fashion 
as follows.  Let $E^0_{X,Y} \in \B_X \vee \B_Y$ be the set
$$ E^0_{X,Y} := \{ (x,y) \in X \times Y: f_1(x,y) \geq \tau^{1/10}; \quad \E( f_2(x,y)^2 | \B_X \vee \B_Y ) \leq \tau \}$$
and define $E^0_{Y,Z} \in \B_Y \vee \B_Z$ and $E^0_{Z,X} \in \B_Z \vee \B_X$ similarly.  We first observe that $f$ is small outside of $E^0_{X,Y}$.  Indeed we have (by the $\B_X \vee \B_Y$-measurability of $E^0_{X,Y}$)
\begin{align*}
\int_X \int_Y f(x,y) 1_{(E^0_{X,Y})^c}(x,y)\ d\mu_X(x) d\mu_Y(y)
&= \int_X \int_Y f_1(x,y) 1_{(E^0_{X,Y})^c}(x,y)\ d\mu_X(x) d\mu_Y(y)\\
&\leq \int_{f_1(x,y) < \tau^{1/10}} f_1(x,y)\ d\mu_X(x) d\mu_Y(y) \\
&\quad\quad+ \int_{\E( f_2(x,y)^2 | \B_X \vee \B_Y ) > \tau}\ d\mu(X) d\mu(Y) \\
&\leq \tau^{1/10} + \frac{1}{\tau} \int_X \int_Y \E( f_2(x,y)^2 | \B_X \vee \B_Y )\ d\mu(X) d\mu(Y) \\
&= \tau^{1/10} + \frac{1}{\tau} \|f_2\|_{L^2(X \times Y)}^2 \\
&= o_{\tau \to 0}(1).
\end{align*}
Let $A,B,C$ be atoms in $\B_X,\B_Y,\B_Z$ respectively such that $A \times B \subset E^0_{X,Y}$, $B \times C \subset E^0_{Y,Z}$, and $C \times A \subset E^0_{Z,X}$, and consider the local quantity
$$
\Lambda_3((f_1+f_2) 1_{A \times B}, (g_1+g_2) 1_{B \times C}, (h_1+h_2) 1_{C \times A} ).$$
We can estimate this as the sum of a main term
$$ \Lambda_3(f_1 1_{A \times B}, g_1 1_{B \times C}, h_1 1_{C \times A} )$$
and three error terms
$$ O( \Lambda_3( |f_2| 1_{A \times B}, 1_{B \times C}, 1_{C \times A} ) ) +
O( \Lambda_3( 1_{A \times B}, |g_2| 1_{B \times C}, 1_{C \times A} ) ) + O( \Lambda_3( 1_{A \times B}, 1_{B \times C}, |h_2| 1_{C \times A} ) ).$$
By definition of $E^0_{X,Y}, E^0_{Y,Z}, E^0_{Z,X}$, we have $f_1,g_1,h_1 \geq \tau^{1/10}$ on $A \times B, B \times C, C \times A$ respectively, and hence the main term is at least
$$ \tau^{3/10} \Lambda_3(1_{A \times B}, 1_{B \times C}, 1_{C \times A} ).$$
On the other hand, we have by construction
$$ \E( f_2(x,y)^2 | A \times B) \leq \tau$$
and hence by Cauchy-Schwarz
$$ \Lambda_3( |f_2| 1_{A \times B}, 1_{B \times C}, 1_{C \times A} ) \leq \tau^{1/2} \Lambda_3( 1_{A \times B}, 1_{B \times C}, 1_{C \times A} ).$$
Similarly for $g_2$ and $h_2$.  Thus the error terms are $O( \tau^{2/10} )$ of the main term.  If $\tau \ll 1$ is chosen sufficiently small, we thus have the local estimate
$$ \Lambda_3(f_1 1_{A \times B}, g_1 1_{B \times C}, h_1 1_{C \times A} ) = O( \Lambda_3((f_1+f_2) 1_{A \times B}, (g_1+g_2) 1_{B \times C}, (h_1+h_2) 1_{C \times A} );$$
summing this over all $A,B,C$ and using \eqref{fgh} and the positivity of $f_1+f_2, g_1+g_2, h_1+h_2$ we conclude that
$$ 
 \Lambda_3(f_1 1_{E^0_{X,Y}} , g_1 1_{E^0_{Y,Z}} , h_1 1_{E^0_{Z,X}} ) \leq O(\eps) + O( 1 / F(M) )
$$
(compare this with \eqref{fgh-bad}).  Since $f_1,g_1,h_1$ are bounded from below by $\tau^{1/10}$ on these sets, we thus have
$$ \Lambda_3(1_{E^0_{X,Y}}, 1_{E^0_{Y,Z}}, 1_{E^0_{Z,X}} ) \leq O(\tau^{-3/10} \eps) + O( \tau^{-3/10} / F(M) ).$$
Now let $E_{X,Y}$ be the subset of $E^0_{X,Y}$, defined as the union of all products $A \times B \subset E^0_{X,Y}$ of atoms 
$A \in \B_X$, $B \in \B_Y$ of size at least $\mu_X(A), \mu_Y(B) \geq \tau / 2^M$.  Since $\B_X$ has complexity at most $M$, the union of all atoms in $\B_X$ of measure at most $\tau/2^M$ has measure at most $\tau$, and thus we see that
$$ \mu_X \times \mu_Y( E^0_{X,Y} \backslash E_{X,Y} ) = O(\tau)$$
and hence from preceding computations
$$ \int_X \int_Y f(x,y) 1_{E_{X,Y}^c}(x,y)\ d\mu_X(x) d\mu_Y(y) = o_{\tau \to 0}(1).$$
We define $E_{Y,Z}, E_{Z,X}$ similarly and observe similar bounds.  Now suppose that the expression
$1_{E_{X,Y}}(x,y) 1_{E_{Y,Z}}(y,z) 1_{E_{Z,X}}(z,x)$ does not vanish identically, then there exist atoms $A,B,C$ of $\B_X,\B_Y,\B_Z$
with $A \times B \subset E_{X,Y}$, $B \times C \subset E_{Y,Z}$, and $C \times A \subset E_{Z,X}$.  In particular
$$  \Lambda_3(1_{A \times B}, 1_{B \times C}, 1_{C \times A} ) \leq O(\tau^{-3/10} \eps) + O( \tau^{-3/10} / F(M) ).$$
On the other hand we have
$$ \Lambda_3(1_{A \times B}, 1_{B \times C}, 1_{C \times A} ) = \mu_X(A) \mu_Y(B) \mu_Z(C) \geq (\tau / 2^M)^3.$$
If we define $F(M) := \lfloor 2^{3M} / \tau^3\rfloor + 1$, and assume that $\eps$ is sufficiently large depending on $\tau$ (noting that $M = O_{F,\tau}(1) = O_\tau(1)$), we obtain a contradiction.  Thus we see that $1_{E_{X,Y}}(x,y) 1_{E_{Y,Z}}(y,z) 1_{E_{Z,X}}(z,x)$ vanishes identically whenever $\eps$ is sufficiently small depending on $\tau$.  If we then set $\tau$ to be a sufficiently slowly decaying function of $\eps$, the claim follows.
\end{proof}

Observe that the actual decay rate $o_{\eps \to 0}(1)$ obtained by the above proof is very slow (it decays like the reciprocal of the inverse tower-exponential function).  It is of interest to obtain better bounds here; it is not known what the exact rate
should be, although the Behrend example (Proposition \ref{behrend}) does show that the decay cannot be polynomial in nature.

The above arguments extend (with some nontrivial difficulty) to hypergraphs, and to proving Szemer\'edi's theorem for progressions
of length $k > 3$; the $k=4$ case was handled in \cite{frankl}, \cite{frankl02} (see also \cite{gowers-hyper-4} for a more recent proof), and the general case in \cite{rodl}, \cite{rodl2}, \cite{rs}, \cite{nrs} and \cite{gowers-reg} (see also \cite{tao:hyper}, \cite{tao:correspondence} for more recent proofs).  We sketch the $k=4$ arguments here (broadly following the ideas from \cite{tao:hyper}, \cite{tao:correspondence}).  Finding progressions of length $4$ in a set $A$ is equivalent to solving the simultaneous relations
$$
\begin{array}{lllll}
&-x_2&-2x_3&-3x_4 &\in A\\
x_1& & -x_3&-2x_4 &\in A\\
2x_1&+x_2&&-x_4 &\in A\\
3x_1&+2x_2 &+ x_3 & &\in A.
\end{array}
$$
Because of this, it is not hard to modify the above arguments to deduce the $k=4$ case of Szemer\'edi's theorem from the 
following lemma:

\begin{lemma}[Strong tetrahedron removal lemma, several variable version]\label{tetra} Let $(X_1,\mu_{X_1}), \ldots, (X_4,\mu_{X_4})$ be probability spaces, and for $ijk=123,234,341,412$ let $f_{ijk}: X_i \times X_j \times X_k \to [0,1]$ be measurable functions such that
$$ \Lambda_4(f_{123},f_{234},f_{341},f_{412}) \leq \eps$$
for some $0 < \eps < 1$, where $\Lambda_4$ is the trilinear form
$$ \Lambda_4(f_{123},f_{234},f_{341},f_{412}) := \int_{X_1} \ldots \int_{X_4} \prod_{ijk=123,234,341,412} f_{ijk}(x_i,x_j,x_k)\ d\mu_{X_1}(x_1) \ldots d\mu_{X_4}(x_4).$$
Then for each $ij=12,23,34,41,13,24$ there exists factors $\B_{ij}$ of $X_i \times X_j$ of complexity at most $O_{\eps}(1)$ and sets $E_{ijk} \in \B_{ij} \vee \B_{ik} \vee \B_{jk}$ for $ijk=123,234,341,412$ with
$\prod_{ijk=123,234,341,412} 1_{E_{ijk}}(x_i,x_j,x_k)$ vanishing identically, such that
$$ \int_{X_i} \int_{X_j} \int_{X_k} f_{ijk}(x_i,x_j,x_k) 1_{E_{ijk}^c}(x_i,x_j,x_k)\ d\mu_{X_1}(x_1) d\mu_{X_2}(x_2) d\mu_{X_3}(x_3) \leq o_{\eps \to 0}(1).$$
\end{lemma}

One can recast this lemma as a statement concerning $3$-uniform hypergraphs; see for instance \cite{tao:hyper}.  We will however not pursue this interpretation here (but see \cite{frankl}, \cite{frankl02}, \cite{rodl}, \cite{rodl2}, \cite{rs}, \cite{nrs}, \cite{gowers-reg}, and \cite{gowers-hyper-4} for a treatment of this material from a hypergraph perspective).

In the case of the triangle removal lemma, it was the $\Box^2$ norm which controlled the size of $\Lambda_4$.  Now the role is played
by the $\Box^3$ norm, defined for a measurable bounded function $f(x,y,z): X \times Y \times Z \to \R$ of three variables by the formula
\begin{align*}
\|f\|_{\Box^3(X \times Y \times Z)}^8 := \int_X \int_X &\int_Y \int_Y \int_Z \int_Z f(x,y,z) f(x,y,z') f(x,y',z) f(x,y',z')\\ &f(x',y,z) f(x',y,z') f(x',y',z) f(x',y',z') \ d\mu_X(x) d\mu_X(x') d\mu_Y(y) d\mu_Y(y') d\mu_Z(z) d\mu_Z(z').
\end{align*}
By modifying the previous arguments we see that the $\Box^3$ norm is indeed a norm (after equating functions that agree almost everywhere) and that we have the generalised von Neumann inequality
$$ |\Lambda_4(f,g,h,k)| \leq \min( \|f\|_{\Box^3}, \|g\|_{\Box^3}, \|h\|_{\Box^3}, \|k\|_{\Box^3} ).$$
The analogue of Lemma \ref{corstruct} is

\begin{lemma}[Lack of uniformity implies correlation with structure]\label{corstruct2}  Let $f: X \times Y \times Z \to [-1,1]$ be such that $\|f\|_{\Box^3} \geq \eta$ for some $\eta > 0$.  Then there exists $A_{X,Y} \subset X \times Y$, $A_{Y,Z} \subset Y \times Z$, and $A_{Z,X} \in Z \times X$ 
$$ |\int_X \int_Y \int_Z 1_{A_{X,Y}}(x,y) 1_{A_{Y,Z}}(y,z) 1_{A_{Z,X}}(z,x) f(x,y,z)\ d\mu_X(x) d\mu_Y(y) d\mu_Z(z)| \geq \eta^8/8.$$
\end{lemma}

This ultimately leads to the following regularity lemma:

\begin{lemma}[Simultaneous Szemer\'edi regularity lemma]\label{ssrl2}  For $ijk=123,234,341,412$, let $f_{ijk}: X_i \times X_j \times X_k \to [0,1]$ be measurable, let $\tau > 0$, and let $F: \N \to \N$ be an arbitrary increasing function (possibly depending on $\tau$).  Then there exists an integer $M = O_{F,\tau}(1)$, factors $\B_{ij}$ of $X_i \times X_j$ of complexity at most $M$ for $ij=12,23,34,41,13,24$ and decompositions 
$f_{ijk} = f_{ijk,1} + f_{ijk,2} + f_{ijk,3}$ for $ijk=123,234,341,412$ where
\begin{itemize}
\item ($f_{ijk,1}$ is structured) We have $f_{ijk,1} = \E(f_{ijk}|\B_{ij} \vee \B_{jk} \vee \B_{ik})$.
\item ($f_{ijk,2}$ is small) We have $\|f_{ijk,2}\|_{L^2(X_i \times X_j \times X_k)} \leq \tau$.
\item ($f_{ijk,3}$ is very uniform) We have $\|f_{ijk,3}\|_{\Box^3(X_i \times X_j \times X_k)} \leq 1/F(M)$.
\item (Positivity) $f_{ijk,1}$ and $f_{ijk,1} + f_{ijk,2}$ take values in $[0,1]$.
\end{itemize}
\end{lemma}

One would then like to repeat the proof of Lemma \ref{trl4} by applying this lemma to decompose each function $f_{ijk}$ into three components $f_{ijk,1}$, $f_{ijk,2}$, $f_{ijk,3}$, and then somehow eliminate the latter two terms to reduce to the structured component $f_{ijk,1}$.  The reason for doing this is that, as $f_{ijk,1}$ is measurable with respect to the bounded complexity factor $\B_{ij} \vee \B_{jk} \vee \B_{ik}$, one can decompose this function (which is a function of three variables $x_i,x_j,x_k$) as a polynomial combination of functions of just two variables (or more precisely, as a linear combination of functions of the form $f_{ij}(x_i,x_j) f_{jk}(x_j,x_k) f_{ik}(x_i,x_k)$).  One can then apply a (slight generalisation of) the triangle removal lemma to
handle such functions; more generally, the strategy is to deduce these sort of removal lemmas for functions of $k$ variables, from similar lemmas concerning functions of $k-1$ variables.  In executing this strategy, there is little difficulty in disposing
of the very uniform components $f_{ijk,3}$, if one takes advantage of the freedom to make the growth function $F$ extremely rapid (one needs to take $F$ to be tower-exponential or faster, to counteract the very weak decay present in the two-variable removal
lemmas).  To dispose of the small components $f_{ijk,2}$ takes a little more work, however.  In the above arguments, one implicitly
used the independence of the underlying factors $\B_X, \B_Y, \B_Z$.  In the current situation, the factors $\B_{ij}$ are not independent of each other, which makes it difficult to eliminate the $f_{ijk,2}$ factors directly.  However, this can be addressed by applying the (two-variable) regularity lemma to simultaneously regularise all the atoms in the factors $\B_{ij}$, making them essentially indepenent relative to one-variable factors.  As one might imagine, making this strategy rigorous is somewhat delicate, and in particular the various large and small parameters (such as $\tau$ and $F$) that appear in the regularity lemmas need to be chosen correctly.
See for instance \cite{tao:hyper} for one such realisation of this type of argument.  More recently, an infinitary approach, using a correspondence principle similar in spirit to the Furstenberg correspondence principle, has been employed to give a slightly different proof of the above results, in which the various large and small parameters in the argument have been set to infinity or zero, thus leading to a cleaner (but less elementary) version of the argument; see \cite{tao:correspondence}.  

\section{Relative triangle removal}

The triangle removal result proven in the previous section, Lemma \ref{trl}, only has non-trivial content when the underlying graph $G$ is \emph{dense}, or more precisely when it contains more than $o_{\eps \to 0}(n^2)$ edges, since otherwise one could simply delete all the edges in $G$ to remove the triangles.  This is related to the fact that Lemma \ref{trl} only implies the existence of progressions of length three in dense sets of integers, but not in sparse sets.  However, it is a remarkable and useful fact that
results such as Lemma \ref{trl}, which ostensibly only apply to dense objects, can in fact be extended ``for free'' to sparse objects, as long as the sparse object has large \emph{relative} density with respect to a sufficiently \emph{pseudorandom} object.  This type of ``transference principle'' from the dense category to the relatively dense category was the decisive new ingredient in the result in \cite{gt-primes} that the primes contained arbitrarily long arithmetic progressions.  We will not prove that result here, however we present a simplified version of that result which already captures many of the key ideas.

If $n \geq 1$ is an integer and $0 \leq p \leq 1$, let $G(n,p)$ be the standard Erd\H{o}s-Renyi random graph on $n$ vertices $\{1,\ldots,n\}$, in which each pair of vertices defines an edge in $G(n,p)$ with an identical independent probability of $p$.

\begin{proposition}[Relative triangle removal lemma]\label{rtrl}\cite{klr}, \cite{tao:gauss}  Let $n > 1$ and $1/\log n \leq p \leq 1$, let $0 < \eps < 1$, and let $H = G(n,p)$.  Then with probability $1 - o_{n \to \infty;\eps}(1)$ the following claim is true: whenever $G$ is a subgraph of $H$ which contains fewer than $\eps p^3 n^3$ triangles, then it is possible to delete $o_{\eps \to 0}(p^2 n^2) + o_{n \to \infty;\eps}(p^2 n^2)$ edges from $G$ to create a new graph $G$ which contains no triangles whatsoever.
\end{proposition}

This result in fact extends to much sparser graphs $G(n,p)$, indeed one can take $p = n^{-1/2+\delta}$ for any fixed $0 < \delta < 1/2$; see \cite{klr}.  This argument proceeded by a careful generalisation of the usual regularity lemma to the setting of sparse subsets of pseudorandom graphs.  As one corollary of their result, one can conclude that if $A$ is a random subset of the positive integers with $\P( n \in A ) = n^{-1/2+\delta}$, and with the events $n \in A$ being independent, then almost surely every subset of $A$ of positive density would contain infinitely many progressions of length three.  We shall proceed differently, using a ``soft'' transference argument, inspired by the ergodic theory approach, which follows closely the treatment in \cite{gt-primes} (and also 
\cite{tao:gauss}).  So far, this argument can only handle logarithmic sparsities rather than polynomial, but requires much less randomness on the graph $G(n,p)$; indeed a suitably ``pseudorandom'' graph would also suffice for this argument.  (For the precise definition of the pseudorandomness needed, see \cite{tao:gauss}.)

Let $(X,\mu_X)=(Y,\mu_Y)=(Z,\mu_Z)$ be the vertex set $\{1,\ldots,n\}$ with the uniform distribution.  Fix the random graph $H = G(n,p)$, and let $\nu(x,y)$ be the function on $\{1,\ldots,n\} \times \{1,\ldots,n\}$ which equals $1/p$ when $(x,y)$ lies in $H$ and $0$ otherwise; we can think of $\nu$ as a function on $X \times Y$, $Y \times Z$, or $Z \times X$.  Note from Chernoff's inequality that even though $\nu$ is not bounded by $O(1)$, with probability $1 - o_{n \to \infty}(1)$, $\nu$ has average 
close to $1$:
$$ \int_X \int_Y \nu(x,y)\ d\mu_X(x) d\mu_Y(y) = 1 + o_{n \to \infty}(1).$$
More sophisticated computations of this sort show that many other correlations of $\nu$ with itself are close to $1$.  For instance, one can show that with probability $1 - o_{n \to \infty}(1)$, we have the octahedral correlation estimate
\begin{equation}\label{octo}
\begin{split}
&\int_X \int_X \int_Y \int_Y \int_Z \int_Z \nu(x,y) \nu(x,y') \nu(x',y) \nu(x',y') \\
& \quad \nu(y,z) \nu(y,z') \nu(y',z) \nu(y',z') \\
& \quad \nu(z,x) \nu(z,x') \nu(z',x) \nu(z',x') \\
&\quad\quad\ d\mu_X(x) d\mu_X(x') d\mu_Y(y) d\mu_Y(y') d\mu_Z(z) d\mu_Z(z') = 1 + o_{n \to \infty}(1).
\end{split}
\end{equation}
(In \cite{tao:gauss}, this estimate, together with some simpler versions, are referred to as the \emph{linear forms condition} on $\nu$.)  To prove Proposition \ref{rtrl}, it then suffices to prove the following variant of Lemma \ref{trl4}:

\begin{lemma}[Relative strong triangle removal lemma, several variable version]\label{trl-rel} Let $(X,\mu_X)$, $(Y,\mu_Y)$, $(Z,\mu_Z)$, $\nu$ be as above, and let $0 < \eps \leq 1$.  With probability $1 - o_{n \to \infty;\eps}(1)$, the following claim is true: whenever $f: X \times Y \to [0,1]$, $g: Y \times Z \to [0,1]$, and $h: Z \times X \to [0,1]$ be measurable functions such that $\Lambda_3(f \nu,g \nu,h \nu) \leq \eps$, then there exists factors $\B_X, \B_Y, \B_Z$ of $X,Y,Z$ respectively of complexity at most $O_{\eps}(1)$ and sets $E_{X,Y} \in \B_X \vee \B_Y$, $E_{Y,Z} \in \B_Y \vee \B_Z$, $E_{Z,X} \in \B_Z \vee \B_X$ respectively
with $1_{E_{X,Y}}(x,y) 1_{E_{Y,Z}}(y,z) 1_{E_{Z,X}}(z,x)$ vanishing identically, such that
\begin{align*}
\int_X \int_Y f(x,y) \nu(x,y) 1_{E_{X,Y}^c}(x,y)\ d\mu_X(x) d\mu_Y(y),&\\
\int_Y \int_Z g(y,z) \nu(y,z) 1_{E_{Y,Z}^c}(y,z)\ d\mu_Y(y) d\mu_Z(z),&\\
\int_Z \int_X h(z,x) \nu(z,x) 1_{E_{Z,X}^c}(z,x)\ d\mu_Z(z) d\mu_X(x) &\leq o_{\eps \to 0}(1).
\end{align*}
\end{lemma}

We leave the deduction of Proposition \ref{rtrl} from Lemma \ref{trl-rel} as an exercise.  Note that the only new feature here is the presence of the weight $\nu$, which causes functions such as $f\nu$ to be unbounded.  Nevertheless, it turns out to be possible to use arguments similar to those in the preceding section and obtain this result with a little effort from its unweighted counterpart, Lemma \ref{trl4}.

The first thing to do is to check that the generalised von Neumann inequality, \eqref{gvni}, continues to hold in the weighted setting:

\begin{lemma}[Relative generalised von Neumann inequality]\label{rgvn}\cite{tao:gauss}  Let the notation be as above.  Then with probability $1 - o_{n \to \infty}(1)$, the following claim is true: whenever $f: X \times Y \to \R$, $g: Y \times Z \to \R$ and $h: Z \times X \to \R$ bounded in magnitude by $\nu+1$ (thus for instance $|f(x,y)| \leq \nu(x,y)+1$ for all $(x,y) \in X \times Y$,
then
$$ |\Lambda_3(f, g, h)| \leq 4\min( \|f \|_{\Box^2}, \|g \|_{\Box^2}, \|h \|_{\Box^2} ) + o_{n \to \infty}(1).$$
\end{lemma}

See also \cite{gt-primes} for a closely related computation.  We also remark that the estimate \eqref{gcz} also continues to hold in this setting because that estimate did not require $f$ to be bounded.

\begin{proof}(Sketch only) By symmetry it suffices to show that
$$ |\int_X \int_Y \int_Z f(x,y) g(y,z) h(z,x)\ d\mu_X(x) d\mu_Y(y) d\mu_Z(z)| \leq \|f \|_{\Box^2} + o_{n \to \infty}(1).$$
Note that it is easy to verify that $\|\nu+1\|_{\Box^2} = 2 + o_{n \to \infty}(1)$ with high probability, and hence $\|f\|_{\Box^2} =O(1)$.  We eliminate the $h$ function by Cauchy-Schwarz in the $z,x$ variables and reduce to showing
\begin{align*}
|\int_X \int_Y \int_Y \int_Z f(x,y) f(x,y') g(y,z) g(y',z) (\nu(z,x)+1)\ &d\mu_X(x) d\mu_Y(y) d\mu_Y(y') d\mu_Z(z)| \\
&\leq 8\|f \|_{\Box^2}^2 + o_{n \to \infty}(1)
\end{align*}
and then eliminate $g$ by a Cauchy-Schwarz in the $y,y',z$ variables and reduce to showing
\begin{align*}
|\int_X \int_X \int_Y \int_Y f(x,y) f(x,y') f(x',y) f(x',y') W(x,x',y,y')\ &d\mu_X(x) d\mu_X(x') d\mu_Y(y) d\mu_Y(y')| \\
&\leq 16\|f \|_{\Box^2}^4 + o_{n \to \infty}(1)
\end{align*}
where
$$ W(x,x',y,y') := \int_Z (\nu(y,z)+1) (\nu(y',z)+1) (\nu(z,x)+1) (\nu(z,x')+1)\ d\mu_Z(z).$$
If $W \equiv 16$ then we would be done by definition of the $\Box^2$ norm.  So it suffices to show that
\begin{align*}|\int_X \int_X \int_Y \int_Y &f(x,y) f(x,y') f(x',y) f(x',y') |W(x,x',y,y')-16| \\
&d\mu_X(x) d\mu_X(x') d\mu_Y(y) d\mu_Y(y')| \leq o_{n \to \infty}(1).
\end{align*}
By one last Cauchy-Schwarz this follows from the estimate
\begin{align*}
|\int_X \int_X \int_Y \int_Y &(\nu(x,y)+1) (\nu(x,y')+1) (\nu(x',y)+1) (\nu(x',y')+1) |W(x,x',y,y')-16|^2 \\
&d\mu_X(x) d\mu_X(x') d\mu_Y(y) d\mu_Y(y')| \leq o_{n \to \infty}(1)
\end{align*}
which can be easily verified from correlation estimates such as \eqref{octo}.
\end{proof}

In light of this lemma, we can continue to neglect errors which are small in $\Box^2$ norm as being negligible.  The key to establishing Lemma \ref{trl-rel} now rests with the following decomposition:

\begin{theorem}[Structure theorem]\label{rkoop}\cite{tao:gauss} Let the notation be as above, let $f: X \times Y \to [0,1]$ be a function, and let $\sigma > 0$.  Then there exists a decomposition
$$ f\nu = f_1 + f_2 + f_3$$
where $f_1$ is non-negative and obeys the uniform upper bound 
$$f_1(x,y) \leq 1 \hbox{ for all } (x,y) \in X \times Y,$$
$f_2$ is non-negative and obeys the smallness bound 
\begin{equation}\label{fsig}
\int_X \int_Y f_2(x,y)\ d\mu_X(x) d\mu_Y(y) = o_{n \to \infty;\sigma}(1),
\end{equation}
and $f_3$ obeys the uniformity estimate
\begin{equation}\label{f3}
\|f_3\|_{\Box^2(X \times Y)} = o_{\sigma \to 0}(1).
\end{equation}
Furthermore $f_1+f_3$ is also non-negative.
\end{theorem}

This theorem should be compared with Lemma \ref{srl}.  The key point is that it approximates the function $f\nu$, for which we have no good uniform bounds, for the function $f_1$, which is bounded by $1$.  With this theorem (and Lemma \ref{rgvn}) it is now a simple
matter to deduce Lemma \ref{trl-rel} from Lemma \ref{trl4}:

\begin{proof}[Proof of Lemma \ref{trl-rel}]  We may assume that $n$ is sufficiently large depending on $\eps$, as the claim is trivial otherwise.  Let $0 < \sigma \leq \eps$ be chosen later.
We apply Theorem \ref{rkoop} to decompose $f\nu=f_1+f_2+f_3$, $g \nu=g_1+g_2+g_3$, $h \nu=h_1+h_2+h_3$, thus
$$ \Lambda_3(f_1+f_2+f_3,g_1+g_2+g_3,h_1+h_2+h_3) \leq \eps.$$
Since $f_1+f_3, g_1+g_3, h_1+h_3, f_2, g_2, h_2$ are all non-negative, we conclude
$$ \Lambda_3(f_1+f_3, g_1+g_3, h_1+h_3) \leq \eps.$$
Repeated application of Lemma \ref{rgvn} and \eqref{f3} (and the hypothesis $\sigma \leq \eps$) then gives
$$ \Lambda_3(f_1, g_1, h_1) \leq o_{\eps \to 0}(1).$$
The functions $f_1,g_1,h_1$ are bounded, so we may apply Lemma \ref{trl4} and obtain $\B_X, \B_Y, \B_Z$ of $X,Y,Z$ respectively of complexity at most $O_{\eps}(1)$ and sets $E_{X,Y} \in \B_X \vee \B_Y$, $E_{Y,Z} \in \B_Y \vee \B_Z$, $E_{Z,X} \in \B_Z \vee \B_X$ respectively with $1_{E_{X,Y}}(x,y) 1_{E_{Y,Z}}(y,z) 1_{E_{Z,X}}(z,x)$ vanishing identically, such that
\begin{align*}
\int_X \int_Y f_1(x,y) 1_{E_{X,Y}^c}(x,y)\ d\mu_X(x) d\mu_Y(y),&\\
\int_Y \int_Z g_1(y,z) 1_{E_{Y,Z}^c}(y,z)\ d\mu_Y(y) d\mu_Z(z),&\\
\int_Z \int_X h_1(z,x) 1_{E_{Z,X}^c}(z,x)\ d\mu_Z(z) d\mu_X(x) &\leq o_{\eps \to 0}(1).
\end{align*}
From \eqref{fsig} we have similar estimates for $f_2$, $g_2$, $h_2$:
\begin{align*}
\int_X \int_Y f_2(x,y) 1_{E_{X,Y}^c}(x,y)\ d\mu_X(x) d\mu_Y(y),&\\
\int_Y \int_Z g_2(y,z) 1_{E_{Y,Z}^c}(y,z)\ d\mu_Y(y) d\mu_Z(z),&\\
\int_Z \int_X h_2(z,x) 1_{E_{Z,X}^c}(z,x)\ d\mu_Z(z) d\mu_X(x) &\leq o_{n \to \infty;\sigma}(1).
\end{align*}
Also, from \eqref{f3}, \eqref{gcz} and the complexity bounds on $\B_X, \B_Y, \B_Z$ we have similar estimates for $f_3$, $g_3, h_3$:
\begin{align*}
\int_X \int_Y f_3(x,y) 1_{E_{X,Y}^c}(x,y)\ d\mu_X(x) d\mu_Y(y),&\\
\int_Y \int_Z g_3(y,z) 1_{E_{Y,Z}^c}(y,z)\ d\mu_Y(y) d\mu_Z(z),&\\
\int_Z \int_X h_3(z,x) 1_{E_{Z,X}^c}(z,x)\ d\mu_Z(z) d\mu_X(x) &\leq o_{\sigma \to 0; \eps}(1).
\end{align*}
If we choose $\sigma$ sufficiently small depending on $\eps$, we thus have
\begin{align*}
\int_X \int_Y f(x,y) 1_{E_{X,Y}^c}(x,y)\ d\mu_X(x) d\mu_Y(y),&\\
\int_Y \int_Z g(y,z) 1_{E_{Y,Z}^c}(y,z)\ d\mu_Y(y) d\mu_Z(z),&\\
\int_Z \int_X h(z,x) 1_{E_{Z,X}^c}(z,x)\ d\mu_Z(z) d\mu_X(x) &\leq o_{\eps \to 0}(1) + o_{n \to \infty;\eps}(1)
\end{align*}
and the claim follows.
\end{proof}

Notice how the complexity estimates on $\B_X, \B_Y, \B_Z$ were essential in allowing one to transfer the unweighted triangle removal lemma, Lemma \ref{trl4}, to the weighted setting, Lemma \ref{trl-rel}.

It remains to prove the structure theorem, Theorem \ref{rkoop}.  A full proof (in much greater generality) of this theorem can
be found in \cite{tao:gauss}, while a closely related theorem appears in \cite{gt-primes}.  We give only a brief summary of the
argument here.  Broadly speaking, we follow the energy increment strategy as used to prove Corollary \ref{kv}.  However, we cannot use Lemma \ref{corstruct} as it only applies for functions $f$ which are bounded.  We must therefore redefine the notion of ``structure'', replacing the notion of a tensor product $1_A(x) 1_B(y)$ with the notion of a \emph{dual function} ${\mathcal D}f(x,y)$ of a function $f: X \times Y \to \R$, defined as
$$ {\mathcal D} f(x,y) := \int_X \int_Y f(x,y') f(x',y) f(x',y')\ d\mu_X(x') d\mu_Y(y').$$
Observe that we have the identity
$$ \int_X \int_Y f(x,y) {\mathcal D} f(x,y)\ d\mu_X(x) d\mu_Y(y) = \| f\|_{\Box^2(X \times Y)}^4.$$
Thus if a function $f$ has large $\Box^2$ norm then it correlates with its own dual function.  This fact will be used as a substitute
for Lemma \ref{corstruct}.  One key property of dual functions are that they can be bounded even when $f$ is unbounded; in particular, with high probability we have ${\mathcal D}(\nu+1)$ bounded pointwise by $O(1)$, and hence ${\mathcal D} f$ will also be bounded for any $f$ bounded pointwise in magnitude by $\nu+1$.  Each of these dual functions can define finite factors
$\B_{{\mathcal D} f, \eps}$ for any resolution $\eps > 0$ by partitioning the range of ${\mathcal D} f$ into intervals of length $\eps$ and letting $\B_{{\mathcal D} f, \eps}$ be the factor generated by the inverse image of these intervals.  (For technical reasons it is convenient to randomly shift this partition in order to negate certain boundary effects - which ultimately lead to the small error $f_2$ appearing in Theorem \ref{rkoop} - but let us gloss over this minor detail here.)  Define a \emph{dual factor of complexity $M$ and resolution $\eps$} to be a factor of the form $\B = \B_{{\mathcal D} f_1,\eps} \vee \ldots \vee \B_{{\mathcal D} f_M,\eps}$ where $f_1,\ldots,f_M$ are bounded in magnitude by $\nu+1$.  These factors are the counterparts of the factors $\B_X \vee \B_Y$ studied in the previous section.  A crucial feature of these factors is (with high probability) that the random weight function $\nu$ is uniformly distributed with respect all to these factors; more precisely, with probability $1 - o_{n \to \infty;\eps,M}(1)$ we have $\E(\nu|\B) = 1 + o_{n \to \infty; M,\eps}(1)$ outside of
an exceptional set $\Omega = \Omega_\B$ with $\int_X \int_Y 1_\Omega(x,y) (\nu(x,y)+1)\ d\mu_X d\mu_Y = o_{n \to \infty; M,\eps}(1)$ for all dual factors of complexity $M$.  This fact is somewhat nontrivial to prove; one needs to invoke the Weierstrass approximation theorem to approximate the indicator function of atoms in $\B$ by polynomial combinations of the dual functions ${\mathcal D} f$ (with the approximation being uniform outside of a small exceptional set $\Omega$), and then using tools such as the Gowers-Cauchy-Schwarz inequality  one can control the inner product of $\nu$ with such polynomials.  See \cite{tao:gauss}, \cite{gt-primes} for details.

Once one has these dual factors with respect to which $\nu$ is (essentially) uniformly distributed, one can then develop a counterpart of Lemma \ref{fei}, which roughly speaking asserts that if $f$ is a function bounded in magnitude by $\nu$, and $\B$ is a dual factor of some complexity $M$ and resolution $\eps$ for which $\|f-\E(f|\B)\|_{\Box^2} \geq \eta$, then with high probability one can find an extension $\B'$ of $\B$ which is a dual factor of complexity $M+1$ and resolution $\eps$, for which the energy
$\| \E(f|\B') \|_{L^2}^2$ has increased from $\| \E(f|\B) \|_{L^2}^2$ by some factor $c(\eta) - o_{\eps \to 0}(1) - o_{n \to \infty; M,\eps}(1)$ for some $c(\eta) > 0$.  This is essentially proven by the same Pythagoras theorem argument used to establish Lemma \ref{fei}, though one has to take some care because $f$, being bounded by $\nu$, does not enjoy good $L^2$ bounds (though the conditional expectations $\E(f|\B)$, $\E(f|\B')$ enjoy uniform bounds outside of a small exceptional set).  One can then iterate this as in the proof of Corollary \ref{kv} to obtain Theorem \ref{rkoop} (with some additional $o_{n \to \infty;\eps}(1)$ errors arising from exceptional sets etc. that can be placed in the small error $f_2$).  See \cite{tao:gauss}, \cite{gt-primes} for details.

\section{Szemer\'edi's original proof}\label{szsc}

We now discuss some of the ideas behind Szemer\'edi's original proof \cite{szemeredi} of his theorem.  This is a remarkably subtle combinatorial argument, and there is no chance that we can describe the full argument here, but we can at least begin to motivate part of the argument.  Rather than plunge directly into the full setup of the argument, we will begin with some naive first attempts at the problem, which do not fully work, but which indicate the steps that need to be taken to obtain a full proof.

The task is, given $k \geq 3$, to show that any subset $A$ of integers whose upper density $\delta = \delta[A] := \limsup_{N\to \infty} \frac{|A\cap [-N,N]|}{2N+1}$ is positive contains at least one progression of length $k$.
The first idea dates back to the original argument of Roth \cite{roth} for the $k=3$ case, which is to try to induct downwards on the upper density  of the set (this is known as the \emph{density increment method}).  If $\delta$ is extremely large, say $\delta > 1 - 1/2k$, then the result is easy, because even a randomly chosen progression will have a good chance of being entirely contained in $A$.  Now one assumes inductively that $A$ has some given upper density $\delta > 0$, and that the theorem has already been proven for higher values of $\delta$.  It is not hard to show that the set of $\delta$ for which Szemer\'edi's theorem holds must be open, so if we can verify in this ``maximal bad density'' case\footnote{This trick is vaguely reminiscent of the reduction to minimal topological dynamical systems, or to ergodic measure-preserving systems.  Unfortunately these tricks seem to be mutually exclusive; if one takes sequences of maximal density then it becomes difficult to convert the argument into a dynamical setting.} that progressions of length $k$ exist, then we are done.

Suppose for contradiction that the set $A$ of this critical density $\delta$ did not have any progressions of length $k$, even though all sets of higher density did have progressions.  What this means is that $A$ cannot contain within it arbitrarily large progressions on which $A$ has higher density.  In other words, we cannot find a sequence of progressions $P_1, P_2,\ldots$ in $\Z$ with length tending to infinity for which $\limsup_{n \to \infty} |A \cap P_n|/|P_n| >\delta$, since if this were the case it would not be difficult to piece together out of the $A \cap P_n$ a set with slightly higher upper density than $A$, but which still had no progressions, contradicting the hypothesis on $\delta$.  Thus we must have $\limsup_{n \to \infty} |A \cap P_n|/|P_n| \leq \delta$ whenever $|P_n| \to \infty$.  In other words, we have the upper bound
\begin{equation}\label{ac}
 |A \cap P| \leq (\delta + o_{|P| \to \infty;A}(1)) |P|
\end{equation}
for all progressions $P$.
[Incidentally, if we knew Szemer\'edi's theorem in the first place, one would deduce immediately that the only such sets $A$ are those sets with density $\delta=0$ or density $\delta=1$, but of course we cannot use Szemer\'edi's theorem to prove itself in such a circular manner!]

Thus on a long progression $P$, the density of $A$ cannot significantly exceed $\delta$.  It is still possible for the density of $A$ to be significantly \emph{less} than $\delta$ on such progressions - but this cannot happen too often, as this would (in conjunction with the upper bound) eventually cause $A$ itself to have density less than $\delta$.  This idea can be easily quantified, and leads to the statement is that given any length $N$, the set
$$ \{ n \in \Z: |A \cap [n, n+N)| = (\delta + o_{N \to \infty;A}(1)) N \}$$
has upper density $1 - o_{N\to \infty;A}(1)$.  Thus ``most'' progressions of length $N$ have density $\delta + o_{N \to \infty;A}(1)$.

This then leads to the next idea, which is to partition the integers into \emph{blocks} $[nN,(n+1)N)$ - progressions of length $N$, in which $n$ is a multiple of $N$.  Call such a block \emph{saturated} if it has the expected density $\delta + o_{N \to \infty;A}(1)$, thus most blocks (in an upper density sense) are saturated.  Suppose temporarily that we could in fact assume that \emph{all} blocks are saturated.  Then we could conclude the argument as follows.  We can colour the $n^{th}$ block $[nN,(n+1)N)$ in one of $2^N$ colours depending on how $A$ is situated inside that block; more precisely, we can color the block $[nN,(n+1)N)$ by the set $\{ 0 \leq i < N: nN + i \in A \}$.  Actually we only need $2^N-1$ colours because the block, being saturated, cannot be completely devoid of elements of $A$.  We have thus coloured all the integers into finitely many colours, and hence by van der Waerden's theorem there is a monochromatic progression of blocks of length $k$.  These blocks have $A$ contained in them in identical fashions, and the blocks are not completely devoid of elements of $A$, so it is not hard to see that the progression of blocks induces a progression of elements of $A$ of the same length, and we are done.

Unfortunately, life is not so simple, and we have the unsaturated blocks to deal with.  While the (lower) density of these exceptional blocks is somewhat small in an absolute sense - it is $o_{N \to \infty;A}(1)$ - it is not very small when compared against the number of colours, $2^N-1$ (or against the reciprocal of this number, to be precise).  Van der Waerden's theorem is nowhere near robust enough to handle such a severe influx of ``uncoloured'' elements.  (It can however deal with a rather easy degenerate case in which the density of saturated blocks unexpectedly happens to be incredibly close to $1$, say at least $1-c(N)$ for some explicit but extremely small $c(N) > 0$ whose exact value depends on the constants arising from van der Waerden's theorem.)  Here we encounter a recurring problem in this field: we are always dealing with quantities which are small, but not small enough.   One is always seeking ways to
somehow iteratively improve the smallness, or at least convert the smallness to another type of smallness which is more robust, in
order to get around this basic issue.

Let's try something else for now.  Suppose we can locate $k$ large blocks of integers, say $[0,N), [N, 2N), \ldots, [(k-1)N, kN)$, which are all saturated.  (This is not hard since the upper density of saturated blocks easily exceeds $1-1/2k$ when $N$ is large enough.)  Let's try to find progressions of length $k$ in $A$ with one element in each block.  Suppose we have somehow (presumably by some sort of an inductive hypothesis) managed to already find many progressions of length $k-1$ in $A$ with one element in each of the first $k-1$ of these blocks.  We can extend each of these progressions by one element, which will most likely lie in the final block $[(k-1)N,kN)$.  (Some of them will not.  However observe that $A$ has to be more or less uniformly distributed on any saturated block, because on any sub-interval of proportional size, $A$ has to have density not much larger $\delta$, and thus on subtraction it must have density not much less than $\delta$ either.  Because of this it is very plausible that a significant fraction of the progressions of $k-1$ located from the induction step will have $k^{th}$ element in the final block as claimed.)  Let $B$ denote the set of all such additional elements of these progressions in $[(k-1)N,kN)$.  If we had a lot of progressions of length $k-1$, it is plausible to expect (by simple counting heuristics) that $B$ should have some positive density in $[(k-1)N,kN)$ (indeed, one expects the density to be comparable to $\delta^{k-1}$).  If $B$ intersects $A$, then we are done.

Unfortunately, $B$ and $A$ are both rather sparse sets inside $[(k-1)N,kN)$ - one has density about $\delta^{k-1}$ (assuming some appropriate induction hypothesis), and the other has density about $\delta$.  These are too sparse to force an intersection unconditionally.  However, we do know that $A$ obeys some good uniform distribution bounds on progressions - its density is always bounded from above, and often bounded from below.  This would be useful if $B$ was somehow made out of progressions (or even better, if the \emph{complement} of $B$ was made out of progressions, since upper bounds on the density of $A$ in the complement of $B$ translate to lower bounds on the density of $A$ in $B$), but we do not have such good structural control on $B$ and it could well be just a generic sparse subset of $[(k-1)N,kN)$, and we are stuck. Indeed, there is nothing right now that stops $B$ from simply being some subset of the complement of $A$, and no matter how structured or uniformly distributed $A$ is, we cannot prevent such an event from happening.

Szemer\'edi's ingenious solution to this problem is to \emph{extend} this sequence of $k$ blocks in an additional direction, which gives $B$ (and more importantly, the complement of $B$) enough of an ``arithmetic progression'' structure that one can eventually get lower bounds on the density of $A$ in $B$.  

To get a preliminary idea of how this idea works, suppose that we have a moderately long progression of saturated 
blocks $P_1, \ldots, P_L$, thus we have $P_i = [a + ir, a+ir + N)$ for some $a \in \Z$ and $r \geq N$, and
\begin{equation}\label{acapp}
|A \cap P_i| = (\delta + o_{N \to \infty;A}(1)) N \hbox{ for all } 1 \leq i \leq L.
\end{equation}
Here $L$ is a moderately large number, though it will be smaller than the length $N$ of each block: $1 \leq L \leq N$.  (Given that the set of saturated blocks has upper density $1 - o_{N \to \infty;A}(1)$, it would be unreasonable to hope to obtain a progression of saturated blocks of length comparable to $N$ or more.)  Let us define $A_i \subset [0,N)$ to be the set $A \cap P_i$, translated
backwards by $a+ir$.  

Now let $B \subset [0,N)$ be a set of some size $\alpha N$.  Then heuristically we expect $A_i \cap B$ to have size $\approx \delta \alpha N$. Now, as discussed before, any \emph{individual} $A_i$ need not have any intersection with $B$.  However, once one considers the sequence $A_1,\ldots,A_L$ there is a kind of ``mixing'' phenomenon that forces at least one of the $A_i$ to have at least the right number of elements inside $B$:

\begin{lemma}[Single lower mixing]\label{slm} Let $P_1,\ldots,P_L$ be a progression of saturated blocks, with attendant sets $A_1,\ldots,A_L \subset [0,N)$ and let $B \subset [0,N)$ be a set of cardinality $\alpha N$.  Then there exists $1 \leq i \leq L$ such that
$$ |A_i \cap B| \geq (\alpha \delta - o_{L \to \infty;A}(1) ) N.$$
\end{lemma}

\begin{proof} By summing \eqref{acapp} for $1 \leq i \leq L$ we have
$$ |A \cap \bigcup_{i=1}^L P_i| = (\delta + o_{N \to \infty;A}(1)) NL.$$
On the other hand, the set $\bigcup_{i=1}^L (P_i \backslash (B+a+ir))$ can be viewed as the union of $(1-\alpha)N$ arithmetic progressions of length $L$.  Applying \eqref{ac} on each such progression and taking unions, we obtain
$$ |A \cap \bigcup_{i=1}^L (P_i \backslash (B+a+ir))| \leq (\delta + o_{L \to \infty;A}(1)) (1-\alpha) NL.$$
Subtracting the latter estimate from the former, we obtain
$$ |A \cap \bigcup_{i=1}^L (B+a+ir)| \geq (\delta\alpha - o_{L \to \infty;A}(1) - o_{N \to \infty;A}(1)) NL.$$
Since $L \leq N$, the latter error term can be absorbed into the former. The claim then follows from the pigeonhole principle, noting that $A \cap (B+a+ir)$ is just a translate of $A_i \cap B$.
\end{proof}

We can amplify this result substantially.  Firstly, we may work with multiple sets $B_1,\ldots,B_m$ instead of a single set $B$.

\begin{lemma}[Multiple lower mixing] 
Let $P_1,\ldots,P_L$ be a progression of saturated blocks, with attendant sets $A_1,\ldots,A_L \subset [0,N)$ and let $B_1,\ldots,B_m \subset [0,N)$ be sets of cardinality $\alpha_1 N,\ldots,\alpha_m N$ respectively.  Then there exists $1 \leq i \leq L$ such that
$$ |A_i \cap B_j| \geq (\alpha_j \delta - o_{L \to \infty;A,m}(1) ) N \hbox{ for all } 1 \leq j \leq m.$$
\end{lemma}

\begin{proof} Suppose that this claim failed.  Then for each $1 \leq i \leq L$ there exists a $j$ for which
$$ |A_i \cap B_j| < (\alpha_j \delta - o_{L \to \infty;A,m}(1) ) N.$$
This is an $m$-colouring of $\{1,\ldots,L\}$.  By van der Waerden's theorem, $\{1,\ldots,L\}$ must then contain a monochromatic
progression of length $\omega_{L \to \infty;m}(1)$, where $\omega_{L \to \infty;m}(1) = 1 / o_{L \to \infty;m}(1)$ denotes a quantity which goes to infinity as $L \to \infty$ for any fixed $m$.  But then this contradicts Lemma \ref{slm} if the $o()$ constants are chosen properly.
\end{proof}

\begin{corollary}[Multiple mixing]\label{part}
Let $P_1,\ldots,P_L$ be a progression of saturated blocks, with attendant sets $A_1,\ldots,A_L \subset [0,N)$ and let $B_1,\ldots,B_m \subset [0,N)$ be sets of cardinality $\alpha_1 N,\ldots,\alpha_m N$ respectively.  Then there exists $1 \leq i \leq L$ such that
$$ |A_i \cap B_j| = (\alpha_j \delta + o_{L \to \infty;A,m}(1) ) N \hbox{ for all } 1 \leq j \leq m.$$
\end{corollary}

\begin{proof}  Apply the preceding lemma, but with $m$ replaced by $2m$ and with $B_{j+m} := [0,N) \backslash B_j$ for $1 \leq j \leq m$.
\end{proof}

This type of result is useful when $m$ is small compared with $L$.  Since $L$ is in turn small compared to $N$, this means that we can only hope to exploit this mixing property when the number $m$ of sets that we wish to be uniformly distributed with respect to $A$ is small compared with the size $N$ of the block.  At first glance, this will severely limit the usefulness of this mixing property; however, we can use the Szemer\'edi regularity lemma to get around this problem (the key point being that the complexity of the partition created by the regularity lemma - which will be $m$ - does not depend on the number of underlying vertices, which is essentially $N$):

\begin{proposition}[Graph mixing]\label{graphmix} Let $P_1,\ldots,P_L$ be a progression of saturated blocks, with attendant sets $A_1,\ldots,A_L \subset [0,N)$ and let $G_1,\ldots,G_m \subset [0,N) \times [0,N)$ be bipartite graphs connecting two copies of $[0,N)$.
Then there exists $1 \leq i \leq L$ such that
$$ \sum_{b \in [0,N)} \bigl| |\{ a \in A_i: (a,b) \in G_j \}| - \delta |\{ a \in [0,N): (a,b) \in G_j \} | \bigr|
= o_{L \to \infty;A,m}(N^2) \hbox{ for all } 1 \leq j \leq m.$$
\end{proposition}

This is a remarkably strong assertion that the set $A_i$ becomes uniformly distributed with density $\delta$ on the interval $[0,N)$ for many values of $i$.  Note that the error term is completely uniform in the graphs $G_1,\ldots,G_m$ (although it does depend of course on the number $m$ of graphs involved) and also is independent of $N$ (after normalising out the natural $1/N^2$ factor).

\begin{proof}(Sketch)  By van der Waerden's theorem as before we can reduce to the case $m=1$.  Pick an $\eps > 0$ and apply the Szemer\'edi regularity lemma to $G$ to obtain an $\eps$-regular approximation to $G_1$ induced by a partition of complexity $O_\eps(1)$.  Apply Corollary \ref{part} to estimate the contribution of the approximation to obtain a net error of
$o_{L \to \infty;A,\eps}(N^2) + o_{\eps \to 0}(N^2)$.  The claim then follows by choosing $\eps$ to be a sufficiently slowly decaying function of $L$.  
(One could also proceed here using a weaker regularity lemma such as Corollary \ref{kv}.)
\end{proof}

Let us now informally discuss how one can exploit such strong mixing properties to extend progressions of length $k-1$ to progressions of length $k$.  (Actually, for technical inductive reasons we will also need to extend progressions of length $i-1$ to progressions of length $i$ for $1 \leq i \leq k$; we shall return to this point later.)  Suppose we have a sequence of $k$-tuples
$(P_{1,i}, P_{2,i}, \ldots, P_{k,i})$ of saturated blocks for $1 \leq i \leq L$, where each $k$-tuple is in progression, and furthermore the final blocks $P_{k,1}, \ldots, P_{k,L}$ of each $k$-tuple are also in progression.  We can then define sets $A_{j,i} \subset [0,N)$ for $1 \leq j \leq k$ and $1 \leq i \leq L$ as before by intersecting $A$ with $P_{j,i}$ and then translating back to $[0,N)$.  We also make the assumption that $A$ ``looks the same'' in the non-final blocks $P_{1,i},\ldots,P_{k-1,i}$, in the sense that for any $1 \leq j \leq k-1$, the sets $A_{j,i}$ are in fact independent of $i$.  Suppose also that in each $k$-tuple $(P_{1,i}, P_{2,i}, \ldots, P_{k,i})$, we have found ``many'' ($\gg \delta^{k-1} N^2$, in fact) progressions of length $k$, with the $j^{th}$ element of the progression in $P_{j,i}$, and with the first $k-1$ elements in $A$.  Note that in fact once a single $k$-tuple, say $(P_{1,1}, P_{2,1}, \ldots, P_{k,1})$ has this property, then all $k$-tuples do,
since this property depends only on the distribution of $A$ in the non-final blocks $P_{1,i},\ldots,P_{k-1,i}$ and we are assuming that this distribution is independent of $i$.  Later we shall address the rather important question of \emph{how} one could construct such a strange sequence of $k$-tuples; for now, let us simply assume that such a sequence exists.  This sequence shows that $A$ has many progressions of length $k-1$.  We now show that some of these progressions of length $k-1$ can be extended
to progressions of length $k$ in $A$; this is a model of the key inductive step in Szemer\'edi's argument.  

Consider the sets $A_{1,i},\ldots,A_{k-1,i}, A_{k,i}$ in $[0,N)$, which describe the distribution of $A$ in the $k$-tuple
$(P_{1,i},\ldots,P_{k,i})$.  The first $k-1$ of these sets are independent of $i$, while the final set $A_{k,i}$ varies in $i$; however, because the blocks $P_{k,1},\ldots,P_{k,L}$ the final set $A_{k,i}$ obeys the strong mixing properties described earlier.
By hypothesis, we have many progressions of length $k$ in $[0,N)$, with the $j^{th}$ element of such progressions lying in $A_{j,i}$ for $1 \leq j \leq k-1$.  The $k^{th}$ elements of such progressions can be collected into a subset of $[0,N)$ which we shall call $B$; we can then get a reasonable lower bound on the density of $B$ in $[0,N)$ (roughly speaking, we have $|B| \gg \delta^{k-1} N$).  The objective is to get $B$ to intersect $A_{k,i}$ for at least one $i$, as this will generate a progression of length $k$ in $A$.
But this happens for at least one $i$ if $L$ is large enough (depending on $\delta$, but not on $N$), thanks to Lemma \ref{slm}.
(Note that we did not use the strongest mixing properties available; we will utilise those later.)
Indeed the intersection of $B$ with $A_{k,i}$ will be rather large, and by arguing slightly more carefully one can then show that
the $i^{th}$ $k$-tuple $(A_{1,i},\ldots,A_{k,i})$ will contain quite a large number of progressions of length $k$ ($\gg \delta^k N^2$, in fact).

To summarise, by using the mixing properties, we can convert a long sequence of $k$-tuples of blocks, each of which contain many progressions of length $k-1$ in $A$, into a single $k$-tuple of blocks, which contains many progressions of length $k$ in $A$, provided that we have the following two additional properties:
\begin{itemize}
\item The distribution of $A$ in the $k-1$ non-final blocks of the $k$-tuples is fixed as one moves along the sequence.
\item The final block of the $k$-tuples are in progression as one moves along the sequence.
\end{itemize}

This looks like a promising induction-type step.  However it cannot by itself be iterated to generate progressions of length $k$ unconditionally for two reasons.  Firstly, there is the minor objection that we will need a generalisation of the above statement in which progressions of length $k-1$ and $k$ in $A$ are replaced by progressions of length $i-1$ and $i$ in $A$ for various $1 \leq i \leq k$.  This is not hard to address.  The more important objection is that we will need a way of generating not only individual $k$-tuples of blocks that contain progressions of length (say) $k-1$ in $A$, but entire \emph{sequences} of such $k$-tuples which obey additional structural properties.  

The key to obtaining this type of superstructure atop a $k$-tuple of blocks in \cite{szemeredi} is by passing to a ``coarser'' level, and viewing each block as a single element of $\Z$; the saturated blocks (as well as a subset of the saturated blocks which are known as the ``perfect'' blocks) then become subsets of $\Z$.  These sets in turn have upper densities, and one can also define notions of saturated blocks of these sets, which are thus ``blocks of blocks''.  The point is that the task of finding sequences of $k$-tuples of blocks simplifies, on moving to this coarser scale, to the task of finding sequences of $k$-term progressions, which is easier and in fact will follow once one has a suitable $k$-tuple of saturated blocks at this coarse scale.

The details are very technical, but let us just mention some brief highlights here.  Write $A_0 = A$.  One picks a large number $N_0$ for which
there are lots of saturated blocks of length $N_0$ (the upper density of such blocks should be $1 - o_{N_0 \to \infty;A}(1)$).
We subdivide the integers into blocks of length $N_0$, and identify the set of such blocks again with $\Z$, creating a ``coarse scale'' view of the set $A_0$.  (Objects in the coarse scale will be subscripted by $1$, while objects in the fine scale subscripted by $0$.)  The saturated blocks then form a subset $S_1$ of $\Z$ of upper density close to $1$.  Each element of $S_1$ corresponds to a saturated block, with respect to which $A_0$ is distributed in one of $2^{N_0}$ ways.  This can be viewed as a colouring of $S_1$ into $2^{N_0}$ colours.  One of the colour classes must be somewhat prevalent (in particular, occuring with positive upper density);
we designate this as the ``perfect'' colour, and let $A_1 \subset S_1$ be the associated colour class.  (The precise definition of ``prevalent'' is slightly technical - it is sort of an upper density ``relative'' to $S_1$ - and we omit it here.)  $A_1$ has some upper density $\delta_1$; it is possible (after some notational trickery) to run a density increment argument for $A_1$ and reduce to the case where $A_1$ obeys an analogue of the bound \eqref{ac}.  In particular we can pick a large number $N_1$ (much larger than $N_0$) and construct many saturated blocks of $A_1$ of length $N_1$.  The definition of ``saturated'' is a little technical; we require that these blocks not only contain $A_1$ to approximately the right density (i.e. $\delta_1 + o_{N_1 \to \infty;A_1}(1)$), but also contains $S_1$ to approximately the right density ($1 - o_{N_0 \to \infty;A}(1)$, if $N_1$ is large enough).  This can be done by tinkering with the notion of upper density appropriately, as mentioned briefly before; we omit the details.

Now suppose one has a $k$-tuple $P_1,\ldots,P_k$ of saturated blocks of $A_1$, and suppose that one can find many $k$-term progressions with the $i^{th}$ term in $P_i$ for $1 \leq i \leq k$, and also in $A_1$ for $1 \leq i \leq k-1$.  Specifically, let us suppose that for almost all (e.g. with density $1 - o_{N_0 \to \infty;A}(1)$) of the integers $n$ in the middle third of the
final block $P_i$, that there are many $k$-term progressions ending in $n$ with the first $k-1$ terms in $P_1 \cap A_1$, $P_2 \cap A_1, \ldots, P_{k-1} \cap A_1$ respectively.  Most of these integers $n$ are going to also lie in $S_1$ (since $S_1$ fills almost all of $P_k$), and so there should be no difficulty obtaining an arithmetic progression of such $n$ of some moderate length $L_0$
(which can be a slowly growing function of $N_0$), thus each element of this progression is the final element of a $k$-term progression which is mostly in $A_1$.  Now recall that each integer in this coarse representation corresponds to a block of length $N_0$ in the original fine-scale representation.  Thus this arithmetic progression can be identified with a sequence of $L_0$ $k$-tuples of such blocks, where the final block in each $k$-tuple is in arithmetic progression, and all the other blocks have the ``perfect'' colour.  This is essentially the very structure we need in order to run our inductive step and convert the progressions with $k-1$ elements in $A_0$, to progressions with $k$ elements in $A_0$.

To summarise, by coarsening the scale it is possible to convert $k$-tuples of blocks to sequences of $k$-tuples of blocks (and more generally to a type of ``homogeneous, well-arranged'' family of $k$-tuples, as defined in \cite{szemeredi}).  These sequences can then be traded in via the mixing properties to upgrade short progressions in a set $A$ to longer progressions.  By alternating these two arguments in a moderately sophisticated induction argument (passing from fine scales to coarse scales approximately $2^k$ times), one can start with progressions with $0$ elements in one of the $A$ sets and eventually upgrade to progressions with $k$
elements in the original set $A$.  There are some technical issues at intermediate stages of the argument, when descending a scale in a case when only the first $i$ elements of a progression are guaranteed to have the perfect colour, when it becomes important
that the remaining elements are unsaturated.  To achieve this, the graph mixing properties in Proposition \ref{graphmix} become
essential; the progressions are reinterpreted as edges connecting the elements of one block to another.  We omit the details.

\providecommand{\MR}{\relax\ifhmode\unskip\space\fi MR }
\providecommand{\MRhref}[2]{%
  \href{http://www.ams.org/mathscinet-getitem?mr=#1}{#2}
}
\providecommand{\href}[2]{#2}

 \end{document}